\documentclass{cocv}

\usepackage{graphicx}
\usepackage{graphicx,amssymb,amsmath}
\usepackage{epsfig}
\usepackage{subfigure}
\usepackage{cite}
\usepackage{booktabs}
\usepackage{array}
\usepackage{color}
\usepackage{enumerate}
\usepackage{algorithm}
\usepackage{algorithmic}
\usepackage{float}
\usepackage{multirow}
\usepackage{rotating}
\usepackage[misc]{ifsym}
\usepackage{bm}
\usepackage{mathrsfs}
\usepackage[dvipdfm,pdfstartview=FitH,CJKbookmarks=true,bookmarksnumbered=true,bookmarksopen=true, colorlinks,pdfborder=001,linkcolor=blue,anchorcolor=blue,citecolor=blue]{hyperref}
\usepackage[top=1.2in, bottom=1.1in, left=0.8in, right=0.8in]{geometry}

\newtheorem{assumption}[subsection]{Assumption}
\numberwithin{equation}{section}

\begin{document}
\title{A FE-ADMM algorithm for Lavrentiev-regularized state-constrained elliptic control problem}\thanks{The National Natural Science Foundation of China (11571061)}
\author{Zixuan Chen}\address{School of Mathematical Sciences, Dalian University of Technology, Dalian, Liaoning 116025, China;\\
\email{chenzixuan@mail.dlut.edu.cn\ \&\ songxiaoliang@mail.dlut.edu.cn\ \&\ zhangxp@dlut.edu.cn\ \&\ yubo@dlut.edu.cn}}
\author{Xiaoliang Song}\sameaddress{1}
\author{Xuping Zhang}\sameaddress{1}
\author{Bo Yu}\sameaddress{1}
\date{The dates will be set by the publisher}
\begin{abstract}
In this paper, elliptic control problems with pointwise box constraints on the state is considered, where the corresponding Lagrange multipliers in general only represent regular Borel measure functions. To tackle this difficulty, the Lavrentiev regularization is employed to deal with the state constraints. To numerically discretize the resulted problem, since the weakness of variational discretization in numerical implementation, full piecewise linear finite element discretization is employed. Estimation of the error produced by regularization and discretization is done. The error order of full discretization is not inferior to that of variational discretization because of the Lavrentiev-regularization. Taking the discretization error into account, algorithms of high precision do not make much sense. Utilizing efficient first-order algorithms to solve discretized problems to moderate accuracy is sufficient. Then a heterogeneous alternating direction method of multipliers (hADMM) is proposed. Different from the classical ADMM, our hADMM adopts two different weighted norms in two subproblems respectively. Additionally, to get more accurate solution, a two-phase strategy is presented, in which the primal-dual active set (PDAS) method is used as a postprocessor of the hADMM. Numerical results not only verify error estimates but also show the efficiency of the hADMM and the two-phase strategy.
\end{abstract}
%
%
\maketitle
\section{Introduction}
\label{intro}
In this paper, we consider the following elliptic PDE-constrained optimal control problem with box constraints on the state
\begin{equation}\label{eqn:original problem}
           \qquad \left\{ \begin{aligned}
        &\min \limits_{(y,u)\in Y\times U}^{}\ \ J(y,u)=\frac{1}{2}\|y-y_d\|_{L^2(\mathrm{\Omega})}^{2}+\frac{\alpha}{2}\|u\|_{L^2(\mathrm{\Omega})}^{2} \\
        &\qquad{\rm s.t.}\qquad -\mathrm{\mathrm{\Delta}} y=u\ \ \mathrm{in}\  \mathrm{\Omega}, \\
         &\qquad \qquad \qquad y=0\quad  \mathrm{on}\ \mathrm{\Gamma},\\
         &\qquad \qquad\qquad a\leq y(x)\leq b \ \ {\rm a. e.  }\  \mathrm{on}\ \mathrm{\Omega},
                          \end{aligned} \right.\tag{$\mathrm{P}$}
\end{equation}
where $Y:=H_0^1(\mathrm{\Omega})$, $U:=L^2(\mathrm{\Omega})$, $\mathrm{\Omega}\subseteq \mathbb{R}^n\ (n=2,3)$ is a convex, open and bounded domain with $C^{1,1}$- or polygonal boundary $\mathrm{\Gamma}$; the desired state $y_d \in L^2(\mathrm{\Omega})$ is given; $a,\ b\in \mathbb{R}$ and $\alpha>0$ are given parameters. Since the constraints in (\ref{eqn:original problem}) denote closed convex set, (\ref{eqn:original problem}) admits unique solution $(y^*,u^*)$. The solution operator $G$ of the elliptic equation in (\ref{eqn:original problem}) mapping $u$ to $y$ is compact. To be more precise, $G=ES$, where $S: u \rightarrow y$ assigns $u\in L^2(\mathrm{\Omega})$ to the weak solution $y\in H_{0}^{1}(\mathrm{\Omega})$ and $E: H_{0}^{1}(\mathrm{\Omega}) \rightarrow L^2(\mathrm{\Omega})$ is the compact embedding operator. We use ${\rm{(\cdot,\cdot)}}$ to denote the inner product in $L^2(\mathrm{\Omega})$ and use $\|\cdot\|$ to denote the corresponding norm. Through this paper, let us suppose the following Slater condition for (\ref{eqn:original problem}) holds.
\begin{assumption}\label{assumption slater point}
There exists a $\widehat{u}\in L^2(\mathrm{\Omega})$ such that\\
\begin{equation*}
a< (S\widehat{u})(x)< b\ \ \ \ \forall x\in \overline{\mathrm{\Omega}}.
\end{equation*}
\end{assumption}

\begin{rmrk}
Our considerations can also carry over to uniformly elliptic operators
\begin{equation*}
\mathscr{A}y=-\sum^{n}_{i,j=1}\partial_{x_j}(a_{ij}y_{x_j})+c_0y,\qquad a_{ij},c_0\in L^{\infty},\ c_0\geq0,\ a_{ij}=a_{ji}
\end{equation*}
and there is a constant $\theta>0$ such that
\begin{equation*}
\sum_{i,j=1}^{n}a_{ij}(x)\xi_i\xi_j\geq\theta\|\xi\|^2\qquad \mathrm{for\ almost\ all}\ \xi\in\mathbb{R}^n.
\end{equation*}
Boundary condition can also expand to
\begin{equation*}
\partial_n y=0\qquad \mathrm{on}\ \mathrm{\Gamma}.
\end{equation*}
\end{rmrk}

Optimal control problems with state constraints and their numerical realization have been studied extensively recently. Since the Lagrange multiplier associated to (\ref{eqn:original problem}) in general only represents a regular Borel measure (see Casas \cite{casas1993boundary} or Alibert and Raymond \cite{alibert1997boundary}) because of the presence of the pointwise state constraints, the complementarity condition in the optimality conditions cannot be written into a pointwise form. Hence, nonsmooth pointwise reformulations, which are needed in semismooth Newton methods, are not possible. To overcome this difficulty, there are two common approaches, Moreau-Yosida regularization and Lavrentiev regularization. Moreau-Yosida regularization \cite{hintermuller2006feasible, hintermuller2006path} is to convert the state constraint into a penalty term. As in \cite{Pearson2014Preconditioners}, the authors showed that a semismooth Newton method applied to the Moreau-Yosida regularization of (\ref{eqn:original problem}) leads to a $3*3$ block saddle point linear system, whose coefficient matrix is symmetric and indefinite. While in our paper, we focus on the Lavrentiev regularization, whose idea is to replace the state constraint by control-state mixed constraint. We can see from Section \ref{sec:4} that only a $2*2$ block saddle point system has to be solved in each iteration by applying our hADMM, which is based on the inherent structure of the problem.

The Lavrentiev regularized problem has the form of a control-constrained elliptic optimal control problem. As we know, since projection has to be carried out to get the control in each iteration in variational discretization \cite{Hinze2005A}, which means mesh refinement for the control, the error order of the control of variational discretization is generally higher than that of full discretization. However, our error analysis indicates that because of the employment of the Lavrentiev-regularization, the error order of the control of full discretization is not inferior to that of variational discretization, which is the most important reason prompting us to use Lavrentiev regularization. The Lavrentiev regularized problem associated to (\ref{eqn:original problem}) is:
\begin{equation}\label{eqn:modified problem lambda}
           \qquad \left\{ \begin{aligned}
        \min \limits_{(y,u)\in Y\times U}^{}\ \ &J(y,u)=\frac{1}{2}\|y-y_d\|_{L^2(\mathrm{\Omega})}^{2}+\frac{\alpha}{2}\|u\|_{L^2(\mathrm{\Omega})}^{2} \\
        {\rm s.t.}\qquad& -\mathrm{\Delta} y=u\ \ \mathrm{in}\ \mathrm{\Omega}, \\
         &y=0\ \ \ \mathrm{on}\ \mathrm{\Gamma},\\
         &a\leq \lambda u+y \leq b\ \ {\rm a. e.  }\  \mathrm{on}\ \mathrm{\Omega},
                          \end{aligned} \right.\tag{$\mathrm{P_{\lambda}}$}
\end{equation}
where $\lambda>0$ denotes the regularization parameter. Since the constraints in (\ref{eqn:modified problem lambda}) denote closed convex set, (\ref{eqn:modified problem lambda}) admits unique solution $(\overline{y}_\lambda,\overline{u}_\lambda)$. In \cite{meyer2006optimal}, the authors prove the convergence of $(\overline{y}_\lambda,\overline{u}_\lambda)\rightarrow(y^*,u^*)$ in $L^2(\mathrm{\Omega})$ for $\lambda\rightarrow0$. Also, they show that the Lagrange multiplier associated to the mixed control-state constraint in (\ref{eqn:modified problem lambda}) is an $L^2$-function for every $\lambda>0$. In addition, \cite{hinze2010variational} proves the weak convergence of the adjoint states in $L^2$ for $\lambda$ tending to zero and the weak-$*$ convergence of the multipliers in $C(\overline{\mathrm{\Omega}})^*$ to their counterparts of problem (\ref{eqn:original problem}) for $\lambda\downarrow0$. Without loss of generality, we assume that $0<\lambda<1$. We know from \cite{prufert2008convergence} that for the error resulted from Lavrentiev-regularization, the following estimate holds
\begin{equation}
\|u^*-\overline u_{\lambda}\|\leq c\sqrt{\lambda},
\end{equation}
where $c$ is a constant independent of $\lambda$. If we introduce an artificial variable $v=y+\lambda u$, (\ref{eqn:modified problem lambda}) can be transformed into a pure control constrained optimal control problem:
\begin{equation}\label{eqn:modified problem v}
           \left\{ \begin{aligned}
        \min \ \ &J(y,v)=\frac{1}{2}\|y-y_d\|_{L^2(\mathrm{\Omega})}^{2}+\frac{\alpha}{2{\lambda}^2}\|v-y\|_{L^2(\mathrm{\Omega})}^{2} \\
        {\rm s.t.}\ \ &-\mathrm{\Delta}{y}+\frac{1}{\lambda}y=\frac{1}{\lambda}v\ \ \mathrm{in}\  \mathrm{\Omega},\\
         &\ y=0\quad  \ \ \mathrm{on}\ \mathrm{\Gamma},\\
         &\ a\leq v \leq b\ \ {\rm a. e.  }\  \mathrm{on}\ \mathrm{\Omega}.
                          \end{aligned} \right.\tag{$\mathrm{\widetilde{P}_{\lambda}}$}
 \end{equation}
Since (\ref{eqn:modified problem v}) is a pure control-constrained problem, it admits a unique Lagrange multiplier in $L^2(\mathrm{\Omega})$ associated to the inequality constraint.

To numerically solve the regularized problems, we use the {\emph{First discretize, then optimize}} approach. With respect to the discrete methods, the variational discretization has been applied in dealing with (\ref{eqn:modified problem lambda}) in \cite{hinze2010variational}, where the authors give the following error estimates.
\begin{equation}\label{v1}
\|u^*-\overline u_{\lambda,h}\|\leq C\left(\sqrt{\lambda}+\frac{1}{\lambda^2}\left(h^2+\frac{1}{\lambda}h^3+\frac{1}{\lambda^2}h^4\right)\right)
\end{equation}
and
\begin{equation}\label{v2}
\|u^*-\overline u_{\lambda,h}\|\leq C\left(\sqrt{\lambda}+\max\{h|\log(h)|,h^{2-\frac{n}{2}}\}\right),
\end{equation}
where $n=2, 3$ denotes the space dimension and $C$ is a positive constant independent of the finite element grid size $h$ and regularization parameter $\lambda$.

Although the variational discretization avoids explicit discretization of the controls, it is not convenient to solve the resulted problem because the control still belongs to function space. In fact, in each iteration of variational discretization, the grid has to be divided again, which costs lots of computations and storage and is not easy to be implemented. In this paper, we use the full discretization method, in which both the state and control are discretized by piecewise linear functions. The remarkable advantage of full discretization is that it can transform the problem into a finite dimensional problem with a good structure, which is convenient to be implemented numerically. More importantly, we extend the results of \cite{hinze2010variational} to the full discretization case, which results in the following two error estimates.
\begin{equation}\label{f1}
\|u^*-\overline u_{\lambda,h}\|\leq C\left(\sqrt{\lambda}+\frac{1}{\sqrt{\lambda}}h+\frac{1}{\lambda^2}\left(h^2+\frac{1}{\lambda}h^3+\frac{1}{\lambda^2}h^4\right)\right)
\end{equation}
and
\begin{equation}\label{f2}
\|u^*-\overline u_{\lambda,h}\|\leq C\left(\sqrt{\lambda}+\max\{h|\log(h)|,h^{2-\frac{n}{2}}\}\right),
\end{equation}
where $n=2, 3$ denotes the space dimension and $C$ is a positive constant independent of $\lambda$ and $h$.

Although at first glance, the precision of (\ref{v1}) is higher than (\ref{f1}) from the view of $h$, actually it depends on the matching relation between $\lambda$ and $h$. For example, we take $h=2^{-9}$, which is small enough in general. Meanwhile, we take $\lambda=10^{-4}$, where $\lambda$ often has to be smaller in practice. In this case, $\frac{h^2}{\lambda^2}=\frac{h}{\lambda^{\frac{3}{2}}}\cdot\frac{h}{\sqrt{\lambda}}$ is bigger than $\frac{h}{\sqrt{\lambda}}$. In addition, the second error estimate (\ref{f2}) is the same as (\ref{v2}). So it does not mean that the error order of full discretization is inferior to that of variational discretization because of the effect of $\lambda$, i.e. the employment of Lavrentiev regularization, especially when $\lambda$ is very small.

An algorithm called the primal-dual active set method (PDAS) has been used in solving the Lavrentiev-regularized state constrained elliptic control problems in \cite{meyer2007two}, which was proved to be a special semismooth Newton method in \cite{hintermuller2002primal}. Benefiting from the local superlinear convergence rate, semismooth Newton method is a prior choice for solving nonsmooth optimization problem. The error of utilizing numerical methods to solve PDE constrained problem consists of two parts: discretization error and the error of algorithm for discretized problem. The error order of piecewise linear finite element method is $O(h)$, so algorithms of high precision do not make much sense because the discretization error account for the main part. Taking the precision of discretization error into account, using fast first-order algorithm is a wise choice. Actually, using algorithms of high precision will not reduce the error but waste computations and storage. In addition, it is seen in Section \ref{sec:4} that in general we have to solve a $4*4$ block equation system in each iteration, which makes the calculation very large, especially when the finite element grid size $h$ is very small. In \cite{Porcelli2016Preconditioning}, the authors give a method to transform the $4*4$ block equation system to a $2*2$ block one, however, it brings additional computation for the inverse of the mass matrix.

As we know, there are many first order algorithms being used to solve finite dimensional large scale optimization fast, such as accelerated proximal gradient (APG) method \cite{jiang2012inexact, beck2009fast, toh2010accelerated, tseng2008accelerated} and alternating direction method of multipliers (ADMM)\cite{li2015qsdpnal, li2016schur, chen2015efficient, fazel2013hankel, boyd2011distributed}. Motivated by the success of these first order algorithms, an APG method in function space (called Fast Inexact Proximal (FIP) method) was proposed to solve the elliptic optimal control problem involving $L^1$-control cost in \cite{schindele2016proximal}. It is known that whether the APG method is efficient depends closely on whether the step-length is close enough to the Lipschitz constant, however, the Lipschitz constant is not easy to estimate in usual. So in this paper, we focus on ADMM, which was originally proposed in \cite{chan1978finite, gabay1976dual} and has been used broadly in many areas. First, we give a brief overview of ADMM for the following linearly constrained convex optimization problem
\begin{equation}\label{eqn:general convex problems1}
        \left\{\begin{aligned}
        \min \quad &\theta_1(x)+\theta_2(y)\\
        {\rm{s.t.}}\quad  &Ax+By=b,\\
        &x\in \mathcal{X},\\
        &y\in \mathcal{Y},
                          \end{aligned} \right.
\end{equation}
where $\theta_1(x):\mathbb{R}^{n_1}\rightarrow \mathbb{R}$ and $\theta_2(y):\mathbb{R}^{n_2}\rightarrow \mathbb{R}$ are convex functions, $A \in \mathbb{R}^{m\times n_1}$, $B \in \mathbb{R}^{m\times n_2}$ and $b \in \mathbb{R}^{m}$, $\mathcal{X}\subset \mathbb{R}^{n_1}$ and $\mathcal{Y}\subset \mathbb{R}^{n_2}$ are given closed, convex sets. The augmented Lagrangian function of (\ref{eqn:general convex problems1}) is
\begin{equation}\label{eqn:augmented Lagrangian function}
  \mathcal{L}_\sigma(x,y,\lambda;\sigma)=\theta_1(x)+\theta_2(y)+(\lambda,Ax+By-b)+\frac{\sigma}{2}\|Ax+By-b\|^2,
\end{equation}
where $\lambda \in \mathbb{R}^{m}$ is the Lagrange multiplier and $\sigma>0$ is a penalty parameter. Each iteration of ADMM has three main steps
\begin{equation}\label{classical ADMM}
        \left\{\begin{aligned}
      &x^{k+1}={\text{argmin}} \left\{\mathcal{L}_\sigma(x,y^k,\lambda^k;\sigma) {\big|}\ x\in \mathcal{X}\right\},\\
      &y^{k+1}={\text{argmin}} \left\{\mathcal{L}_\sigma(x^{k+1},y,\lambda^k;\sigma) {\big|}\ y\in \mathcal{Y}\right\},\\
      &\lambda^{k+1}=\lambda^k+\sigma(Ax^{k+1}+By^{k+1}-b).
                          \end{aligned} \right.
\end{equation}
The advantage of ADMM is that it separates $\theta_1(x)$ and $\theta_2(y)$ into two subproblems, which makes each subproblem in (\ref{classical ADMM}) could be solved easily. The ADMM algorithm for solving (\ref{eqn:general convex problems1}) has global convergence and sublinear convergence rate at least under some general assumptions.

To apply ADMM type algorithm to (\ref{eqn:modified problem lambda}), we introduce an artificial variable $v=\lambda u+y$, which results in
\begin{equation}\label{eqn:modified problem vu}
           \left\{ \begin{aligned}
        \min \ \ &J(y,v)=\frac{1}{2}\|y-y_d\|_{L^2(\mathrm{\Omega})}^{2}+\frac{\alpha}{2}\|u\|_{L^2(\mathrm{\Omega})}^{2} \\
        {\rm s.t.}\ \ &-\mathrm{\Delta} y=u\ \ \mathrm{in}\  \mathrm{\Omega}, \\
         &\ y=0\quad  \mathrm{on}\ \mathrm{\Gamma},\\
         &\ v-\lambda u-y=0\ \ {\rm a. e.  }\  \mathrm{on}\ \mathrm{\Omega},\\
         &\ a\leq v \leq b\ \ {\rm a. e.  }\  \mathrm{on}\ \mathrm{\Omega}.
                          \end{aligned} \right.\tag{$\mathrm{\widehat{P}_{\lambda}}$}
\end{equation}
Compared with (\ref{eqn:modified problem lambda}) and (\ref{eqn:modified problem v}), (\ref{eqn:modified problem vu}) separates the smooth and nonsmooth terms, which makes it more efficiently to take advantage of ADMM.

The ADMM type algorithm has been used in elliptic optimal control problem with control constraints. In \cite{song2016two}, the authors proposed a heterogeneous ADMM (hADMM) algorithm. The hADMM algorithm employs two different weighted norms in the augmented term in two subproblems respectively, which is different from the classical ADMM. Also, the authors proved the global convergence and the iteration complexity results $o(\frac{1}{k})$. Inspired by the simpleness, facility for implementation and global convergence rate of the hADMM, we employ it to fully discretized Lavrentiev-regularized problem. Although Lavrentiev-regularized problem can be transformed into a pure control-constrained problem as form (\ref{eqn:modified problem v}), it will become ill-conditioned when lambda is very small. Thus we do not apply hADMM to (\ref{eqn:modified problem v}), we use its well structure as reference and apply it to (\ref{eqn:modified problem vu}), which possesses well structure as we see in Section \ref{sec:4}. For the first subproblem of hADMM, it is equivalent to solve a $2*2$ block equation system in each iteration, while using PDAS has to solve a $4*4$ block equation system which should be carefully formed based on the active sets in each iteration. For the second subproblem of hADMM, the solution has a closed form, which is very easy to compute.

 Moreover, to satisfy the need for more accurate solution, a two-phase strategy is also presented, in which the primal-dual active set (PDAS) method is used as a postprocessor of the hADMM algorithm. It is shown in Section \ref{sec:5} that to get a solution of the same precision, the hADMM algorithm and the two-phase strategy are obviously faster than PDAS method respectively.

The paper is organized as follows. Full discretization is considered in Section \ref{sec:2}. Section \ref{sec:3} gives the error estimates of the fully discretized Lavrentiev-regularized problem. In Section \ref{sec:4}, we give the frame of the hADMM algorithm and the PDAS method employed to the discretized problems. Two numerical examples are given to verify the error estimates and the efficiency of the proposed algorithm in Section \ref{sec:5}. Section \ref{sec:6} contains a brief summary of this paper.

\section{Full Finite Element Discretization}
\label{sec:2}
In order to tackle (\ref{eqn:modified problem lambda}) and (\ref{eqn:modified problem vu}) numerically, we consider the full discretization, in which both the state $y$ and the control $u$ are discretized by continuous piecewise linear functions, for which we make the following assumptions.
\begin{assumption}\label{assumption finite element}
$\mathrm{\Omega}\subseteq \mathbb{R}^n$ denotes a bounded domain, $\overline{\mathrm{\Omega}}=\bigcup^{nt}_{j=1}\overline{T}_j$ with admissible quasiuniform sequences of partitions $\{T_j\}^{nt}_{j=1}$ of $\mathrm{\Omega}$, i.e. with $h_{nt}:=\max_j\ {\rm{diam}}(T_j)$ and $\sigma_{nt}:=\min_j{\rm{sup\ diam}} (K);$ $K\subseteq T_j$ there holds $c\leq\frac{h_{nt}}{\sigma_{nt}}\leq C$ uniformly in $nt$ with positive constants $0<c\leq C<\infty$ independent of $nt$. We abbreviate $\tau_h:=\{T_j\}^{nt}_{j=1}$ and set $h=h_{nt}$.
Let $\bar{\mathrm{\Omega}}_h=\bigcup_{T\in \tau_h}T$. We use ${\mathrm{\Omega}}_h$ and $\mathrm{\Gamma}_h$ denoting its interior and boundary respectively. In the case that $\mathrm{\Omega}$ is a convex polyhedral domain, there holds $\mathrm{\Omega}=\mathrm{\Omega}_h$. In the case that $\mathrm{\Omega}$ has a $C^{1,1}$- boundary $\mathrm{\Gamma}$, $\bar{\mathrm{\Omega}}_h$ is convex, whose boundary vertices are all contained in $\mathrm{\Gamma}$, such that
\begin{equation*}
  |\mathrm{\Omega}\backslash {\mathrm{\Omega}}_h|\leq k h^2,
\end{equation*}
where $|\cdot|$ denotes the measure of the set and $k>0$ is a constant.
\end{assumption}

The weak formulation of the state equation involved in (\ref{eqn:modified problem lambda}) and (\ref{eqn:modified problem vu})
\begin{equation}\label{eqn:state equation}
\begin{aligned}
   -\mathrm{\Delta} y&=u \qquad \mathrm{in}\  \mathrm{\Omega}, \\
   y&=0 \qquad \mathrm{on}\ \mathrm{\Gamma}
\end{aligned}
\end{equation}
is given by
\begin{equation}\label{weak formulation:state equation}
  (\nabla y,\nabla z)=(u,z),\qquad \forall z \in H_0^1(\mathrm{\Omega}).
\end{equation}
Let a finite dimensional subspace $Z_h$ of $H^1_0(\mathrm{\Omega})$
\begin{equation}
          Z_h=\left\{z_h\in C(\overline{\mathrm{\Omega}})\ |\ {z_h|}_T \in P_1\ \ \ \forall T\in T_h\ and \ z_h=0\ in\ \overline{\mathrm{\Omega}}\backslash \mathrm{\Omega}_h \right\}
\end{equation}
be the discrete space, where $\mathcal{P}_1$ denotes the space of polynomials whose degree are less than or equal to $1$. Let $\{\phi_i(x)\} _{i=1}^{N_h}$ be a basis of $Z_h$ which satisfies the following properties:
\begin{equation}\label{basic functions properties}
  \phi_i(x) \geq 0,\qquad\|\phi_i(x)\|_{\infty} = 1, \quad \forall i=1,2,...,N_h,\qquad\sum\limits_{i=1}^{N_h}\phi_i(x)=1,
\end{equation}
then (\ref{weak formulation:state equation}) implies that the weak formulation is satisfied for all basis functions $\{\phi_i(x)\} _{i=1}^{N_h}$, i.e.
\begin{equation}
          (\nabla y,\nabla \phi_i)=(u,\phi_i),\qquad \forall i=1,\cdots,N_h.
\end{equation}
We discretize $y(x)$ and $u(x)$ by the same basis of $Z_h$, i.e.
\begin{equation}
          y_h(x)=\sum^{N_h}_{i=1} y_i\phi_i(x)\ \mathrm{and}\ u_h(x)=\sum^{N_h}_{i=1} u_i\phi_i(x),
\end{equation}
where $y_h(x_i)=y_i$ and $u_h(x_i)=u_i$. Then the discrete version of problem (\ref{eqn:modified problem lambda}), (\ref{eqn:modified problem v}) and (\ref{eqn:modified problem vu}) are denoted by (\ref{eqn:discretized problem lambda}), (\ref{eqn:discretized problem v}) and (\ref{eqn:discretized problem vu}) respectively,
\begin{equation}\label{eqn:discretized problem lambda}
  \left\{ \begin{aligned}
        \min\ \ &J_h(y_h,u_h)=\frac{1}{2}\|y_h-y_{d}\|_{L^2(\mathrm{\Omega}_h)}^{2}+\frac{\alpha}{2}\|u_h\|_{L^2(\mathrm{\Omega}_h)}^{2}\\
        {\rm{s.t.}}\ \ &(\nabla y_h, \nabla z_h)=(u_h,z_h)\ \ \forall z_h\in Z_h,  \\
          &a\leq \lambda u_h(x)+y_h(x)\leq b\ \ {\rm a. e.  }\  \mathrm{on}\ \mathrm{\Omega},
                          \end{aligned} \right.\tag{$\mathrm{P}_{\lambda,h}$}
\end{equation}
\begin{equation}\label{eqn:discretized problem v}
  \left\{ \begin{aligned}
        \min\ \ &J_h(y_h,v_h)=\frac{1}{2}\|y_h-y_{d}\|_{L^2(\mathrm{\Omega}_h)}^{2}+\frac{\alpha}{2{\lambda}^2}\|v_h-y_h\|_{L^2(\mathrm{\Omega}_h)}^{2}\\
        {\rm{s.t.}}\ \ &(\nabla y_h, \nabla z_h)+\frac{1}{\lambda}(y_h,z_h)=\frac{1}{\lambda}(v_h,z_h)\ \ \forall z_h\in Z_h,  \\
          &a\leq v_h(x) \leq b\ \ {\rm a. e.  }\  \mathrm{on}\ \mathrm{\Omega},
                          \end{aligned} \right.\tag{$\mathrm{\widetilde{P}}_{\lambda,h}$}
\end{equation}
\begin{equation}\label{eqn:discretized problem vu}
  \left\{ \begin{aligned}
        \min\ \ &J_h(y_h,u_h)=\frac{1}{2}\|y_h-y_{d}\|_{L^2(\mathrm{\Omega}_h)}^{2}+\frac{\alpha}{2}\|u_h\|_{L^2(\mathrm{\Omega}_h)}^{2}\\
        {\rm{s.t.}}\ \ &(\nabla y_h, \nabla z_h)=(u_h,z_h)\ \ \forall z_h\in Z_h, \\
        &v_h-\lambda u_h-y_h=0\ \ {\rm a. e.  }\  \mathrm{on}\ \mathrm{\Omega},\\
          &a\leq v_h(x) \leq b\ \ {\rm a. e.  }\  \mathrm{on}\ \mathrm{\Omega}.
                          \end{aligned} \right.\tag{$\mathrm{\widehat{P}}_{\lambda,h}$}
\end{equation}

\section{Error estimates}
\label{sec:3}
In this section, we extend the results of \cite{hinze2010variational}. The essential difference between \cite{hinze2010variational} and the present paper is that the discretization method in \cite{hinze2010variational} is variational discretization while this paper considers full discretization, in which both the state and control are discretized by piecewise linear functions. The greatest difficulty that full discretization introduces to the error analysis is that the solution of continuous problem is not feasible for discretized problem. To tackle with this difficulty, we utilize the quasi-interpolation operator and complete the error analysis. It is well known that since projection has to be carried out to get the control in each iteration in variational discretization, which means mesh refinement for the control, the error order of the control of variational discretization is generally higher than that of full discretization. However, the error analysis in this section indicates that the error order of the control of full discretization is not inferior to that of variational discretization because of the employment of the Lavrentiev-regularization. In this section, we give two different error estimates, the first one of which depends on $\lambda$ while the second one of which is uniform in $\lambda$.
\subsection{Error estimate for fixed $\lambda$}
For the error analysis below, we have to use a quasi-interpolation operator ${\rm\Pi}_h:L^2(\mathrm{\Omega})\rightarrow Z_h$, which is defined by
\begin{equation*}
        {\rm\Pi}_{h}v=\sum^{N_h}_{i=1} \pi_i(v)\phi_i(x),\ \ \pi_i(v)=\frac{\int_{\mathrm{\Omega}_h} v(x)\phi_i(x) dx}{\int_{\mathrm{\Omega}_h} \phi_i(x) dx},\quad \forall v\in L^2(\mathrm{\Omega}).
\end{equation*}
Let
\begin{equation*}
           V_{ad}=\{v\in L^2(\mathrm{\Omega})~|~a\leq v \leq b\ \ {\rm a. e.  }\  \mathrm{on}\ \mathrm{\Omega}\}
\end{equation*}
and
\begin{equation*}
           V_{ad,h}=\{v_h=\sum_{i=1}^{N_h}v_i\phi_i(x)|~a\leq v_i\leq b\ \ {\rm a. e.  }\  \mathrm{on}\ \mathrm{\Omega}\},
\end{equation*}
then there holds
\begin{equation*}
        v\in V_{ad}\Rightarrow {\rm\Pi}_h v\in V_{ad,h},\ \ \forall v\in L^2(\mathrm{\Omega}).
\end{equation*}
For the interpolation error, the following lemma holds, whose proof can be found in \cite{carstensen1999quasi, de2008finite}.
\begin{lmm}\label{interpolation error lemma}
There exists a constant $C$ independent of $h$ such that
\begin{equation*}
h\|v-{\rm\Pi}_h v\|_{L^2}+\|v-{\rm\Pi}_h v\|_{H^{-1}}\leq Ch^2\|v\|_{H^1}\qquad \forall v\in H^1(\mathrm{\Omega}).
\end{equation*}
\end{lmm}

First we consider the following variational equation
\begin{equation}\label{state equation}
(\nabla w, \nabla z)+\frac{1}{\lambda}(w,z)=(g,z),\qquad \forall z\in H^1_0(\mathrm{\Omega})
\end{equation}
and its discrete version:
\begin{equation}\label{discrete state equation}
(\nabla w_h, \nabla z_h)+\frac{1}{\lambda}(w_h,z_h)=(g,z_h),\qquad \forall z_h\in Z_h,
\end{equation}
where $g\in L^2(\mathrm{\Omega})$. We use $w(g)$ and $w_h(g)$ to denote the solution of (\ref{state equation}) and (\ref{discrete state equation}) respectively, then the following lemma holds.

\begin{lmm}\label{wh(g)-w(g)}
Under Assumption \ref{assumption finite element}, there exists a constant $C(\mathrm{\Omega})$ independent of $\lambda$ such that
\begin{equation*}
\|w_h(g)-w(g)\|_{L^2(\mathrm{\Omega})}\leq C(\mathrm{\Omega})\left(h^2+\frac{1}{\lambda}h^3+\frac{1}{\lambda^2}h^4\right)\|w(g)\|_{H^2(\mathrm{\Omega})}
\end{equation*}
holds true.
\end{lmm}
\begin{proof}
Let $z=w_h(g)-I_hw(g)$ in (\ref{state equation}) and $z_h=w_h(g)-I_hw(g)$ in (\ref{discrete state equation}), then we get
\begin{equation*}
\begin{aligned}
&(\nabla z(g), \nabla (z_h(g)-I_hz(g)))+\frac{1}{\lambda}(z(g),z_h(g)-I_hz(g))=(g,z_h(g)-I_hz(g)),\\
&(\nabla z_h(g), \nabla (z_h(g)-I_hz(g)))+\frac{1}{\lambda}(z_h(g),z_h(g)-I_hz(g))=(g,z_h(g)-I_hz(g)),
\end{aligned}
\end{equation*}
where $I_h$ denotes the linear interpolation operator. Subtracting two equalities above, we arrive at
\begin{equation}
(\nabla (z_h(g)-z(g)), \nabla (z_h(g)-I_hz(g)))+\frac{1}{\lambda}(z_h(g)-z(g),z_h(g)-I_hz(g))=0,
\end{equation}
so
\begin{equation*}
\begin{aligned}
&\ \ \|z_h(g)-z(g)\|^2_{H^1(\mathrm{\Omega})}\\
&\leq\left(\nabla (z_h(g)-z(g)), \nabla (z_h(g)-z(g))\right)+\frac{1}{\lambda}\left(z_h(g)-z(g),z_h(g)-z(g)\right)\\
&=\left(\nabla (z_h(g)-z(g)), \nabla (I_hz(g)-z(g))\right)+\frac{1}{\lambda}\left(z_h(g)-z(g),I_hz(g)-z(g)\right)\\
&\leq \frac{1}{2}\|\nabla(z_h(g)-z(g))\|^2+\frac{1}{2}\|\nabla(I_hz(g)-z(g))\|^2+\frac{1}{\lambda}\left(z_h(g)-z(g),I_hz(g)-z(g)\right)\\
&\leq \frac{1}{2}\|z_h(g)-z(g)\|^2_{H^1(\mathrm{\Omega})}+\frac{1}{2}\|I_hz(g)-z(g)\|^2_{H^1(\mathrm{\Omega})}+\frac{1}{\lambda}\left(z_h(g)-z(g),I_hz(g)-z(g)\right),
\end{aligned}
\end{equation*}
where we have used $\frac{1}{\lambda}>1$ .Then we arrive at
\begin{equation*}
\frac{1}{2}\|z_h(g)-z(g)\|^2_{H^1(\mathrm{\Omega})}\leq\frac{1}{2}\|I_hz(g)-z(g)\|^2_{H^1(\mathrm{\Omega})}+\frac{1}{\lambda}(z_h(g)-z(g),I_hz(g)-z(g)),
\end{equation*}
so
\begin{equation}
\begin{aligned}
&\|z_h(g)-z(g)\|^2_{H^1(\mathrm{\Omega})}\leq\|I_hz(g)-z(g)\|^2_{H^1(\mathrm{\Omega})}+\frac{2}{\lambda}(z_h(g)-z(g),I_hz(g)-z(g))\\
&\qquad\qquad\qquad\qquad\ \ \leq\left(\|z(g)-I_hz(g)\|_{H^1(\mathrm{\Omega})}+\frac{1}{\lambda}\|z(g)-I_hz(g)\|\right)^2.
\end{aligned}
\end{equation}
Standard interpolation error estimates imply
\begin{equation}
\begin{aligned}
&\|z_h(g)-z(g)\|_{H^1(\mathrm{\Omega})}\leq\|z(g)-I_hz(g)\|_{H^1(\mathrm{\Omega})}+\frac{1}{\lambda}\|z(g)-I_hz(g)\|\\
&\qquad\qquad\qquad\qquad\ \ \leq C(\mathrm{\Omega})\left(h+\frac{1}{\lambda}h^2\right)\|z(g)\|_{H^2(\mathrm{\Omega})}.
\end{aligned}
\end{equation}
Let $\phi$ be the solution of
\begin{equation}\label{phi definition}
(\nabla \phi,\nabla z)+\frac{1}{\lambda}(\phi,z)=(w-w_h,z),\qquad \forall z\in H^1_0(\mathrm{\Omega})
\end{equation}
and we have
\begin{equation}\label{zero}
(\nabla(w-w_h),\nabla z_h)+\frac{1}{\lambda}(w-w_h,z_h)=0,\qquad \forall z_h\in Z_h.
\end{equation}
Let $z=w-w_h$ in (\ref{phi definition}) and $z_h=I_h\phi$ in (\ref{zero}), we arrive at
\begin{equation*}
\begin{aligned}
&\|w-w_h\|^2=(\nabla\phi,\nabla(w-w_h))+\frac{1}{\lambda}(\phi,w-w_h)-(\nabla I_h\phi,\nabla(w-w_h))-\frac{1}{\lambda}(I_h\phi,w-w_h)\\
&\qquad\qquad\ \ =(\nabla(\phi-I_h\phi),\nabla(w-w_h))+\frac{1}{\lambda}(\phi-I_h\phi,w-w_h)\\
&\qquad\qquad\ \ \leq\|w-w_h\|_{H^1}\cdot\|\phi-I_h\phi\|_{H^1}+\frac{1}{\lambda}\|w-w_h\|_{H^1}\cdot\|\phi-I_h\phi\|\\
&\qquad\qquad\ \ \leq\|w-w_h\|_{H^1}\cdot Ch\|\phi\|_{H^2}+\frac{1}{\lambda}\|w-w_h\|_{H^1}\cdot Ch^2\|\phi\|_{H^2}\\
&\qquad\qquad\ \ \leq\|w-w_h\|_{H^1}\cdot Ch\|w-w_h\|+\frac{1}{\lambda}\|w-w_h\|_{H^1}\cdot Ch^2\|w-w_h\|,
\end{aligned}
\end{equation*}
where we have used the fact that $\|\phi\|_{H^2}\leq C\|w-w_h\|_{L^2}$. 
Then we arrive at
\begin{equation}
\begin{aligned}
\|w-w_h\|&\leq C(h+\frac{1}{\lambda}h^2)\|w-w_h\|_{H^1}\\
&\leq C(\mathrm{\Omega})(h+\frac{1}{\lambda}h^2)\left(h+\frac{1}{\lambda}h^2\right)\|z(g)\|_{H^2(\mathrm{\Omega})}\\
&\leq C(\mathrm{\Omega})\left(h^2+\frac{1}{\lambda}h^3+\frac{1}{\lambda^2}h^4\right)\|z(g)\|_{H^2(\mathrm{\Omega})}.
\end{aligned}
\end{equation}
\end{proof}

Let $(\overline{y}_{\lambda},\overline{v}_{\lambda})$ and $(\overline{y}_{\lambda,h},\overline{v}_{\lambda,h})$ be the solutions of {\rm{(\ref{eqn:modified problem v})}} and {\rm{(\ref{eqn:discretized problem v})}} respectively, then the optimal system of (\ref{eqn:modified problem v}) is:
\begin{subequations}\label{eqn:KKT for PV}
\begin{eqnarray}
       &&(\nabla \overline{y}_{\lambda}, \nabla z)+\frac{1}{\lambda}(\overline{y}_{\lambda},z)=\frac{1}{\lambda}(\overline{v}_{\lambda},z), \qquad  \forall z\in H_0^1(\mathrm{\Omega})\label{eqn1:the optimal system of PV},\\
       &&(\nabla p_{\lambda}, \nabla z)+\frac{1}{\lambda}(p_{\lambda},z)=(\overline{y}_{\lambda}-y_d+\frac{\alpha}{\lambda^2}(\overline{y}_{\lambda}-\overline{v}_{\lambda}),z),\ \ \forall z\in H_0^1(\mathrm{\Omega}),\label{eqn2:the optimal system of PV}\\
       &&\overline{v}_{\lambda}\in V_{ad}, \qquad (\overline{v}_{\lambda}-\overline{y}_{\lambda}+\frac{\lambda}{\alpha}p_{\lambda},v-\overline{v}_{\lambda})\geq0, \qquad  \forall v\in V_{ad},\label{eqn3:the optimal system of PV}
\end{eqnarray}
\end{subequations}
where $p_{\lambda}$ denotes the adjoint state. And the optimal system of (\ref{eqn:discretized problem v}) is:
\begin{subequations}\label{eqn:KKT for PV_h}
\begin{eqnarray}
       &&(\nabla \overline{y}_{\lambda,h}, \nabla z_h)+\frac{1}{\lambda}(\overline{y}_{\lambda,h},z_h)=\frac{1}{\lambda}(\overline{v}_{\lambda,h},z_h), \qquad  \forall z_h\in Z_h,\label{eqn1:the optimal system of PV_h}\\
       &&(\nabla p_{\lambda,h}, \nabla z_h)+\frac{1}{\lambda}(p_{\lambda,h},z_h)=(\overline{y}_{\lambda,h}-y_d+\frac{\alpha}{\lambda^2}(\overline{y}_{\lambda,h}-\overline{v}_{\lambda,h}),z_h),\ \forall z_h\in Z_h,\label{eqn2:the optimal system of PV_h}\\
       &&\overline{v}_{\lambda,h}\in V_{ad,h}, \qquad (\overline{v}_{\lambda,h}-\overline{y}_{\lambda,h}+\frac{\lambda}{\alpha}p_{\lambda,h},v-\overline{v}_{\lambda,h})\geq 0, \qquad  \forall v\in V_{ad,h},\label{eqn3:the optimal system of PV_h}
\end{eqnarray}
\end{subequations}
where $p_{\lambda,h}$ denotes the adjoint state. Additionally, Let $y(v)$, $y_h(v)$, $p(v)$, $p^h(v)$ and $p_h(v)$ be the solution of
\begin{eqnarray}
        &&(\nabla y,\nabla z)+\frac{1}{\lambda}(y,z)=\frac{1}{\lambda}(v,z),\ \ \forall z\in H_0^1(\mathrm{\Omega}),\label{eqn:y(v)}\\
        &&(\nabla y_h,\nabla z_h)+\frac{1}{\lambda}(y_h,z_h)=\frac{1}{\lambda}(v,z_h),\ \ \forall z_h\in Z_h,\label{eqn:y_h(v)}\\
        &&(\nabla p,\nabla z)+\frac{1}{\lambda}(p,z)=(y(v)-y_d+\frac{\alpha}{\lambda ^2}(y(v)-v),z),\ \ \forall z\in H_0^1(\mathrm{\Omega}),\label{eqn:p(v)}\\
        &&(\nabla p^h,\nabla z_h)+\frac{1}{\lambda}(p^h,z_h)=(y(v)-y_d+\frac{\alpha}{\lambda ^2}(y(v)-v),z_h),\ \ \forall z_h\in Z_h,\label{eqn:p^h(v)}\\
        &&(\nabla p_h,\nabla z_h)+\frac{1}{\lambda}(p_h,z_h)=(y_h(v)-y_d+\frac{\alpha}{\lambda ^2}(y_h(v)-v),z_h),\ \ \forall z_h\in Z_h,\label{eqn:p_h(v)}
\end{eqnarray}
respectively, then we have $\overline y_{\lambda}=y(\overline v_{\lambda}),\ p_{\lambda}=p(\overline v_{\lambda}),\ \overline y_{\lambda,h}=y_h(\overline v_{\lambda,h}),\ p_{\lambda,h}=p_h(\overline v_{\lambda,h})$. The following corollary can be easily derived from Lemma \ref{wh(g)-w(g)}.
\begin{crllr}\label{error estimate 1 auxiliary lemma}
Suppose that Assumption \ref{assumption finite element} is fullfilled. Then there exists a constant $C(\mathrm{\Omega})$ independent of $\lambda$ such that the following estimate is valid
\begin{equation*}
\|y_h(\overline v_{\lambda})-\overline y_{\lambda}\|\leq C(\mathrm{\Omega})\left(h^2+\frac{1}{\lambda}h^3+\frac{1}{\lambda^2}h^4\right)
\end{equation*}
In addition
\begin{equation*}
\lambda\|p^h(\overline v_{\lambda})-p_{\lambda}\|\leq C(\alpha,\mathrm{\Omega})\left(h^2+\frac{1}{\lambda}h^3+\frac{1}{\lambda^2}h^4\right)
\end{equation*}
holds true with a constant $C(\alpha,\mathrm{\Omega})$ independent of $\lambda$.
%
\end{crllr}

\begin{thrm}\label{An error estimate for fixed lambda}
Suppose that Assumption \ref{assumption finite element} is fulfilled. Let $(\overline{y}_{\lambda},\overline{v}_{\lambda})$ and $(\overline{y}_{\lambda,h},\overline{v}_{\lambda,h})$ be the solutions of {\rm{(\ref{eqn:modified problem v})}} and {\rm{(\ref{eqn:discretized problem v})}} respectively, then there exists a constant $C(\alpha,\mathrm{\Omega},\lambda_{max})$ independent of $\lambda$ such that
\begin{equation*}
\|\overline{u}_{\lambda}-\overline{u}_{\lambda,h}\|+\|\overline{y}_{\lambda}-\overline y_{\lambda,h}\|_{H^1(\mathrm{\Omega})}\leq C(\alpha,\mathrm{\Omega},\lambda_{max})\left(\frac{1}{\sqrt{\lambda}}h+\frac{1}{\lambda^2}\left(h^2+\frac{1}{\lambda}h^3+\frac{1}{\lambda^2}h^4\right)\right)
\end{equation*}
is satisfied.
\end{thrm}
\begin{proof}
Because the solution $\overline{v}_{\lambda,h}$ of (\ref{eqn:discretized problem v}) is feasible for (\ref{eqn:modified problem v}), we can insert $\overline{v}_{\lambda,h}$ in (\ref{eqn3:the optimal system of PV}), which gives
\begin{equation}\label{variational equality 1}
        (\overline{v}_{\lambda}-\overline{y}_{\lambda}+\frac{\lambda}{\alpha}p_{\lambda},\overline{v}_{\lambda,h}-\overline{v}_{\lambda})\geq0.
\end{equation}
Let $\widetilde{v}={\rm\Pi}_h \overline{v}_{\lambda}$, where ${\rm\Pi}_h$ is the quasi-interpolation operator defined above. Then $\widetilde{v}$ is feasible for (\ref{eqn:discretized problem v}) and we can insert $\widetilde{v}$ in (\ref{eqn3:the optimal system of PV_h}), which gives
\begin{equation}\label{variational equality 2}
        (\overline{v}_{\lambda,h}-\overline{y}_{\lambda,h}+\frac{\lambda}{\alpha}p_{\lambda,h},\widetilde{v}-\overline{v}_{\lambda})+(\overline{v}_{\lambda,h}-\overline{y}_{\lambda,h}+\frac{\lambda}{\alpha}p_{\lambda,h},\overline{v}_{\lambda}-\overline{v}_{\lambda,h})\geq 0.
\end{equation}
Adding (\ref{variational equality 1}) and (\ref{variational equality 2}) then yields
\begin{equation*}
        (\overline{v}_{\lambda,h}-\overline{y}_{\lambda,h}+\frac{\lambda}{\alpha}p_{\lambda,h},\widetilde{v}-\overline{v}_{\lambda})+(\overline{v}_{\lambda}-\overline{v}_{\lambda,h}-(\overline{y}_{\lambda}-\overline{y}_{\lambda,h})+\frac{\lambda}{\alpha}(p_{\lambda}-p_{\lambda,h}),\overline{v}_{\lambda,h}-\overline{v}_{\lambda})\geq 0.
\end{equation*}
We can rewrite the inequality above into the following form:
\begin{equation}\label{eqn:original error estimate}
\begin{aligned}
        0\leq&(\overline{v}_{\lambda,h}-\overline{y}_{\lambda,h}+\frac{\lambda}{\alpha}p_{\lambda,h},\widetilde{v}-\overline{v}_{\lambda})-\parallel\overline{v}_{\lambda}-\overline{v}_{\lambda,h}\parallel^2+(y_h(\overline v_{\lambda})-\overline y_{\lambda},\overline v_{\lambda,h}-\overline v_{\lambda})\\
        &+\frac{\lambda}{\alpha}(p_{\lambda}-p^h(\overline v_{\lambda}),\overline{v}_{\lambda,h}-\overline{v}_{\lambda})+\underbrace{\frac{\lambda}{\alpha}(p^h(\overline v_{\lambda})-p_h(\overline v_{\lambda}),\overline{v}_{\lambda,h}-\overline{v}_{\lambda})}_{I_1}\\
        &+\underbrace{(\overline y_{\lambda,h}-y_h(\overline v_{\lambda}),\overline v_{\lambda,h}-\overline v_{\lambda})+\frac{\lambda}{\alpha}(p_h(\overline v_{\lambda})-p_{\lambda,h},\overline{v}_{\lambda,h}-\overline{v}_{\lambda})}_{I_2}.\\
\end{aligned}
\end{equation}
Let $v=\overline v_{\lambda,h},\ z_h=p_h(\overline v_{\lambda})-p_{\lambda,h}\in Z_h$ and $v=\overline v_{\lambda},\ z_h=p_h(\overline v_{\lambda})-p_{\lambda,h}\in Z_h$ in (\ref{eqn:y_h(v)}) respectively. Subtracting the two resulted equalities we get
\begin{equation*}
\begin{aligned}
        &(\nabla \overline y_{\lambda,h}-\nabla y_h(\overline{v}_{\lambda}),\nabla p_h(\overline v_{\lambda})-\nabla p_{\lambda,h})+\frac{1}{\lambda}(\overline y_{\lambda,h}-y_h(\overline v_{\lambda}),p_h(\overline v_{\lambda})-p_{\lambda,h})\\
        &=\frac{1}{\lambda}(\overline v_{\lambda,h}-\overline v_{\lambda},p_h(\overline v_{\lambda})-p_{\lambda,h}).
\end{aligned}
\end{equation*}
Let $v=\overline v_{\lambda},\ z_h=\overline y_{\lambda,h}-y_h(\overline v_{\lambda})\in Z_h$ and $v=\overline v_{\lambda,h},\ z_h=\overline y_{\lambda,h}-y_h(\overline v_{\lambda})\in Z_h$ in (\ref{eqn:p_h(v)}) respectively. Subtracting the two resulted equalities we arrive at
\begin{equation*}
\begin{aligned}
        &(\nabla p_h(\overline v_{\lambda})-\nabla p_{\lambda,h},\nabla \overline y_{\lambda,h}-\nabla \overline y_h(\overline{v}_{\lambda}))+\frac{1}{\lambda}(p_h(\overline v_{\lambda})-p_{\lambda,h},\overline y_{\lambda,h}-y_h(\overline v_{\lambda}))\\
        &=(y_h(\overline v_{\lambda})-\overline y_{\lambda,h}+\frac{\alpha}{\lambda^2}(y_h(\overline v_{\lambda})-\overline y_{\lambda,h}-\overline v_{\lambda}+\overline v_{\lambda,h}),\overline y_{\lambda,h}-y_h(\overline v_{\lambda})).
\end{aligned}
\end{equation*}
So we have
\begin{equation*}
        \frac{1}{\lambda}(\overline v_{\lambda,h}-\overline v_{\lambda},p_h(\overline v_{\lambda})-p_{\lambda,h})=(y_h(\overline v_{\lambda})-\overline y_{\lambda,h}+\frac{\alpha}{\lambda^2}(y_h(\overline v_{\lambda})-\overline y_{\lambda,h}-\overline v_{\lambda}+\overline v_{\lambda,h}),\overline y_{\lambda,h}-y_h(\overline v_{\lambda})).
\end{equation*}
Then we can rewrite $I_2$ in (\ref{eqn:original error estimate}) as
\begin{equation}
\begin{aligned}\label{eqn:B}
       I_2&=(\overline y_{\lambda}-y_h(\overline v_{\lambda}),\overline v_{\lambda,h}-\overline v_{\lambda})+\frac{\lambda}{\alpha}(p_h(\overline v_{\lambda})-p_{\lambda,h},\overline{v}_{\lambda,h}-\overline{v}_{\lambda})\\
          &=(\overline y_{\lambda,h}-y_h(\overline v_{\lambda}),\overline v_{\lambda,h}-\overline v_{\lambda})\\
          &\ +\frac{\lambda^2}{\alpha}(y_h(\overline v_{\lambda})-\overline y_{\lambda,h}+\frac{\alpha}{\lambda^2}(y_h(\overline v_{\lambda})-\overline y_{\lambda,h}-\overline v_{\lambda}+\overline v_{\lambda,h}),\overline y_{\lambda,h}-y_h(\overline v_{\lambda}))\\
          &=-(1+\frac{\lambda^2}{\alpha})\|y_h(\overline v_{\lambda})-\overline y_{\lambda,h}\|^2+2(\overline y_{\lambda,h}-y_h(\overline v_{\lambda}),\overline v_{\lambda,h}-\overline v_{\lambda}).
\end{aligned}
\end{equation}
Similarly let $v=\overline v_{\lambda}, z_h=\overline y_{\lambda,h}-y_h(\overline v_{\lambda})\in Z_h$ in (\ref{eqn:p^h(v)}) and (\ref{eqn:p_h(v)}) respectively. Subtracting the two resulted equalities we derive
\begin{equation*}
\begin{aligned}
       &(\nabla p^h(\overline v_{\lambda})-\nabla p_h(\overline v_{\lambda}),\nabla \overline y_{\lambda,h}-\nabla y_h(\overline v_{\lambda}))+\frac{1}{\lambda}(p^h(\overline v_{\lambda})-p_h(\overline v_{\lambda}),\overline y_{\lambda,h}-y_h(\overline v_{\lambda}))\\
       &=(y(\overline v_{\lambda})-y_h(\overline v_{\lambda})+\frac{\alpha}{\lambda^2}(y(\overline v_{\lambda})-y_h(\overline v_{\lambda})),\overline y_{\lambda,h}-y_h(\overline v_{\lambda})).
\end{aligned}
\end{equation*}
Let $v=\overline v_{\lambda,h}, z_h=p^h(\overline v_{\lambda})-p_h(\overline v_{\lambda})\in Z_h$ and $v=\overline v_{\lambda}, z_h=p^h(\overline v_{\lambda})-p_h(\overline v_{\lambda})\in Z_h$ in (\ref{eqn:y_h(v)}) respectively. Subtracting the two resulted equalities we have
\begin{equation*}
\begin{aligned}
       &(\nabla \overline y_{\lambda,h}-\nabla y_h(\overline v_{\lambda}),\nabla p^h(\overline v_{\lambda})-\nabla p_h(\overline v_{\lambda}))+\frac{1}{\lambda}(\overline y_{\lambda,h}-y_h(\overline v_{\lambda}),p^h(\overline v_{\lambda})-p_h(\overline v_{\lambda}))\\
       &=\frac{1}{\lambda}(\overline v_{\lambda,h}-\overline v_{\lambda},p^h(\overline v_{\lambda})-p_h(\overline v_{\lambda})).
\end{aligned}
\end{equation*}
So we arrive at
\begin{equation*}
        \frac{1}{\lambda}(\overline v_{\lambda,h}-\overline v_{\lambda},p^h(\overline v_{\lambda})-p_h(\overline v_{\lambda}))=(y(\overline v_{\lambda})-y_h(\overline v_{\lambda})+\frac{\alpha}{\lambda^2}(y(\overline v_{\lambda})-y_h(\overline v_{\lambda})),\overline y_{\lambda,h}-y_h(\overline v_{\lambda})).
\end{equation*}
Then we can rewrite $I_1$ in (\ref{eqn:original error estimate}) as
\begin{equation}\label{eqn:A}
\begin{aligned}
       I_1&=\frac{\lambda}{\alpha}(p^h(\overline v_{\lambda})-p_h(\overline v_{\lambda}),\overline{v}_{\lambda,h}-\overline{v}_{\lambda})\\
           &=\frac{\lambda^2}{\alpha}(\overline y_{\lambda}-y_h(\overline v_{\lambda})+\frac{\alpha}{\lambda^2}(\overline y_{\lambda}-y_h(\overline v_{\lambda})),\overline y_{\lambda,h}-y_h(\overline v_{\lambda}))\\
           &=(1+\frac{\lambda^2}{\alpha})(\overline y_{\lambda}-y_h(\overline v_{\lambda}),\overline y_{\lambda,h}-y_h(\overline v_{\lambda})).
\end{aligned}
\end{equation}
Inserting (\ref{eqn:B}) and (\ref{eqn:A}) into (\ref{eqn:original error estimate}), we get
\begin{equation*}
\begin{aligned}
       0\leq &-\|\overline v_{\lambda}-\overline v_{\lambda,h}\|^2+(y_h(\overline v_{\lambda})-\overline y_{\lambda},\overline v_{\lambda,h}-\overline v_{\lambda})+\frac{\alpha}{\lambda}(p_{\lambda}-p^h(\overline v_{\lambda}),\overline v_{\lambda,h}-\overline v_{\lambda})\\
       &-(1+\frac{\lambda^2}{\alpha})(y_h(\overline v_{\lambda})-\overline y_{\lambda,h},y_h(\overline v_{\lambda})-\overline y_{\lambda,h})+2(\overline y_{\lambda,h}-y_h(\overline v_{\lambda}),\overline v_{\lambda,h}-\overline v_{\lambda})\\
       &-(1+\frac{\lambda^2}{\alpha})(\overline y_{\lambda}-y_h(\overline v_{\lambda}),y_h(\overline v_{\lambda})-\overline y_{\lambda,h})+(\overline v_{\lambda,h}-\overline y_{\lambda,h}+\frac{\lambda}{\alpha}p_{\lambda,h},\widetilde{v}-\overline v_{\lambda})\\
       =&-[\|\overline v_{\lambda}-\overline v_{\lambda,h}\|^2-2(\overline y_{\lambda}-\overline y_{\lambda,h},\overline v_{\lambda}-\overline v_{\lambda,h})+\|\overline y_{\lambda}-\overline y_{\lambda,h}\|^2]-\frac{\lambda^2}{\alpha}\|\overline y_{\lambda}-\overline y_{\lambda,h}\|^2\\
       &+(\overline y_{\lambda}-y_h(\overline v_{\lambda}),\overline v_{\lambda,h}-\overline v_{\lambda})+\frac{\alpha}{\lambda}(p_{\lambda}-p^h(\overline v_{\lambda}),\overline v_{\lambda,h}-\overline v_{\lambda})\\
       &-(1+\frac{\lambda^2}{\alpha})(\overline y_{\lambda}-\overline y_{\lambda,h},y_h(\overline v_{\lambda})-\overline y_{\lambda})+(\overline v_{\lambda,h}-\overline y_{\lambda,h}+\frac{\lambda}{\alpha}p_{\lambda,h},\widetilde{v}-\overline v_{\lambda})\\
       =&-\lambda^2\|\overline u_{\lambda}-\overline u_{\lambda,h}\|^2-\frac{\lambda^2}{\alpha}\|\overline y_{\lambda}-\overline y_{\lambda,h}\|^2\\
       &+(\overline y_{\lambda}-y_h(\overline v_{\lambda}),\lambda(\overline u_{\lambda,h}-\overline u_{\lambda})+(\overline y_{\lambda,h}-\overline y_{\lambda}))+\frac{\lambda}{\alpha}(p_{\lambda}-p^h(\overline v_{\lambda}),\lambda(\overline u_{\lambda,h}-\overline u_{\lambda})\\
       &+(\overline y_{\lambda,h}-\overline y_{\lambda}))-(1+\frac{\lambda^2}{\alpha})(y_h(\overline v_{\lambda})-\overline y_{\lambda},\overline y_{\lambda}-\overline y_{\lambda,h})+\frac{\lambda}{\alpha}(\alpha \overline u_{\lambda,h}+p_{\lambda,h},\widetilde{v}-\overline v_{\lambda}).
\end{aligned}
\end{equation*}
So we derive
\begin{equation*}
\begin{aligned}
       \alpha \|\overline u_{\lambda}-\overline u_{\lambda,h}\|^2+\|\overline y_{\lambda}-\overline y_{\lambda,h}\|^2&\leq \frac{\alpha}{\lambda}(\overline y_{\lambda}-y_h(\overline v),\overline u_{\lambda,h}-\overline u_{\lambda})+(p_{\lambda}-p^h(\overline v_{\lambda}),\overline u_{\lambda,h}-\overline u_{\lambda})\\
       &+\frac{1}{\lambda}(p_{\lambda}-p^h(\overline v_{\lambda}),\overline y_{\lambda,h}-\overline y_{\lambda})+(y_h(\overline v_{\lambda})-\overline y_{\lambda},\overline y_{\lambda,h}-\overline y_{\lambda})\\
       &+\frac{1}{\lambda}(\alpha \overline u_{\lambda,h}+p_{\lambda,h},\widetilde{v}-\overline v_{\lambda}).
\end{aligned}
\end{equation*}
Using Young's inequality we get
\begin{equation*}
\begin{aligned}
       &(\alpha-2k)\|\overline u_{\lambda}-\overline u_{\lambda,h}\|^2+(1-2k)\|\overline y_{\lambda}-\overline y_{\lambda,h}\|^2\\
       &\leq (\frac{\alpha^2}{k\lambda^2}+\frac{1}{k})\|\overline y_{\lambda}-y_h(\overline v_{\lambda})\|^2+(\frac{1}{k\lambda^2}+\frac{1}{k\lambda^4})\lambda^2\|p_{\lambda}-p^h(\overline v_{\lambda})\|^2\\
       &+\frac{1}{\lambda}\|\alpha \overline u_{\lambda,h}+p_{\lambda,h}\|_{H^1(\mathrm{\Omega})}\cdot\|\widetilde{v}-\overline v_{\lambda}\|_{H^{-1}(\mathrm{\Omega})},
\end{aligned}
\end{equation*}
with $k>0$ arbitrary. Then Corollary \ref{error estimate 1 auxiliary lemma} and Lemma \ref{interpolation error lemma} yield
\begin{equation}
\begin{aligned}
      &(\alpha-2k)\|\overline u_{\lambda}-\overline u_{\lambda,h}\|^2+(1-2k)\|\overline y_{\lambda}-\overline y_{\lambda,h}\|^2\\
      &\leq C(\alpha,\mathrm{\Omega},\lambda_{max})\left[\frac{1}{k\lambda^4}\left(h^2+\frac{1}{\lambda}h^3+\frac{1}{\lambda^2}h^4\right)^2+\frac{1}{\lambda}h^2\right].
\end{aligned}
\end{equation}
Let $k=\frac{1}{4}\min(\alpha,1)$ to make $\alpha-2k>0$ and $1-2k>0$, then we arrive at
\begin{equation}\label{error for u}
\|\overline u_{\lambda}-\overline u_{\lambda,h}\|\leq C(\alpha,\mathrm{\Omega},\lambda_{max})\left(\frac{1}{\sqrt{\lambda}}h+\frac{1}{\lambda^2}\left(h^2+\frac{1}{\lambda}h^3+\frac{1}{\lambda^2}h^4\right)\right).
\end{equation}
For $\|\overline y_{\lambda}-\overline y_{\lambda,h}\|_{H^1}$, we have $\forall t>0$,
\begin{equation*}
\begin{aligned}
    \|\overline y_{\lambda}-\overline y_{\lambda,h}\|_{H^1}^2\leq& C\{a(\overline{y}_{\lambda}-\overline{y}_{\lambda,h},\overline{y}_{\lambda}-y_h(\overline{u}_{\lambda}))+a(\overline{y}_{\lambda}-\overline{y}_{\lambda,h},y_h(\overline{u}_{\lambda})-\overline{y}_{\lambda,h})\}\\
    =&C\langle\overline{y}_{\lambda}-y_h(\overline{u}_{\lambda}),\overline{u}_{\lambda}-\overline{u}_{\lambda,h}\rangle_{H^1_0,H^{-1}}\\
    &+C\langle y_h(\overline{u}_{\lambda})-\overline{y}_{\lambda}+\overline{y}_{\lambda}-\overline{y}_{\lambda,h},\overline{u}_{\lambda}-\overline{u}_{\lambda,h}\rangle_{H^1_0,H^{-1}}\\
    \leq& Ct\|\overline y_{\lambda}-y_h(\overline{u}_{\lambda})\|_{H^1}^2+\frac{C}{t}\|\overline{u}_{\lambda}-\overline{u}_{\lambda,h}\|^2+Ct\|\overline y_{\lambda}-y_h(\overline{u}_{\lambda})\|_{H^1}^2\\
    &+\frac{C}{t}\|\overline{u}_{\lambda}-\overline{u}_{\lambda,h}\|^2+Ct\|\overline y_{\lambda}-\overline{y}_{\lambda,h}\|_{H^1}^2+\frac{C}{t}\|\overline{u}_{\lambda}-\overline{u}_{\lambda,h}\|^2\\
    =&2Ct\|\overline y_{\lambda}-y_h(\overline{u}_{\lambda})\|_{H^1}^2+Ct\|\overline y_{\lambda}-\overline{y}_{\lambda,h}\|_{H^1}^2+\frac{3C}{t}\|\overline{u}_{\lambda}-\overline{u}_{\lambda,h}\|^2,
\end{aligned}
\end{equation*}
which implies
\begin{equation}
    (1-Ct)\|\overline y_{\lambda}-\overline y_{\lambda,h}\|_{H^1}^2\leq 2Ct\|\overline y_{\lambda}-y_h(\overline{u}_{\lambda})\|_{H^1}^2+\frac{3C}{t}\|\overline{u}_{\lambda}-\overline{u}_{\lambda,h}\|^2.
\end{equation}
We choose $t=\frac{1}{2C}$ to make $1-Ct>0$, then we derive
\begin{equation}
    \|\overline y_{\lambda}-\overline y_{\lambda,h}\|_{H^1}\leq \widehat{C}\{\|\overline y_{\lambda}-y_h(\overline{u}_{\lambda})\|_{H^1}+\|\overline{u}_{\lambda}-\overline{u}_{\lambda,h}\|\}.
\end{equation}
We know from standard error estimates that $\|\overline y_{\lambda}-y_h(\overline{u}_{\lambda})\|_{H^1}\leq Ch\|\overline{u}_{\lambda}\|$, which together with (\ref{error for u}) implies
\begin{equation}
\|\overline y_{\lambda}-\overline y_{\lambda,h}\|_{H^1}\leq C(\alpha,\mathrm{\Omega},\lambda_{max})\left(h+\frac{1}{\sqrt{\lambda}}h+\frac{1}{\lambda^2}\left(h^2+\frac{1}{\lambda}h^3+\frac{1}{\lambda^2}h^4\right)\right).
\end{equation}
Since $0<\lambda<1$, so $h<\frac{1}{\sqrt{\lambda}}h$. Then the term $h$ can be abandoned from the formula above and we can get the assertion.

\end{proof}

\subsection{Error estimate uniform in $\lambda$}
We now derive an error estimate which does not depend on $\lambda$. Let $(\overline{y}_{\lambda},\overline{u}_{\lambda})$ and $(\overline{y}_{\lambda,h},\overline{u}_{\lambda,h})$ denote the solutions of {\rm{(\ref{eqn:modified problem lambda})}} and {\rm{(\ref{eqn:discretized problem lambda})}} respectively, then the optimal system of (\ref{eqn:modified problem lambda}) is:
\begin{subequations}\label{eqn:KKT for P_lambda}
\begin{eqnarray}
       && (\nabla \overline y_\lambda,\nabla z)=(\overline u_\lambda,z)\ \ \forall z\in H_0^1(\mathrm{\Omega}),\label{eqn1:KKT for P_lambda} \\
       &&(\nabla p_\lambda,\nabla z)=(\overline y_\lambda-y_d-\mu_a+\mu_b,z),\ \forall z\in H_0^1(\mathrm{\Omega}),\label{eqn2:KKT for P_lambda} \\
       &&\alpha \overline u_\lambda+p_\lambda+\lambda \mu_b-\lambda \mu_a=0 \quad {\rm a. e.  \ in} \ \mathrm{\Omega},\label{eqn3:KKT for P_lambda}\\
       &&(\mu_a,a-\lambda \overline u_\lambda-\overline y_\lambda)=(\mu_b,\lambda \overline u_\lambda+\overline y_\lambda-b)=0,\label{eqn4:KKT for P_lambda} \\
       &&\mu_a(x)\geq0,\ \ \mu_b(x)\geq0 \quad {\rm a. e.  \ in} \ \mathrm{\Omega} \label{eqn5:KKT for P_lambda}, \\
       &&a\leq\lambda\overline u_\lambda+\overline y_\lambda\leq b \quad\ {\rm a. e.  \ in} \ \mathrm{\Omega}, \label{eqn6:KKT for P_lambda}
\end{eqnarray}
\end{subequations}
where $p_{\lambda}$ is the adjoint state and $\mu_{a},\ \mu_{b}$ are Lagrange multipliers associated to the regularized pointwise state constraints in {\rm{(\ref{eqn:modified problem lambda})}}. Similarly, the optimal system of (\ref{eqn:discretized problem lambda}) is:
\begin{subequations}\label{eqn:KKT for P_lambda,h}
\begin{eqnarray}
       && (\nabla \overline y_{\lambda,h},\nabla z_h)=(\overline u_{\lambda,h},z_h)\ \ \forall z_h\in Z_h,\label{eqn1:KKT for P_lambda,h} \\
       &&(\nabla p_{\lambda,h},\nabla z_h)=(\overline y_{\lambda,h}-y_d-\mu_{a,h}+\mu_{b,h},z_h),\ \forall z_h\in Z_h,\label{eqn2:KKT for P_lambda,h} \\
       &&\alpha \overline u_{\lambda,h}+p_{\lambda,h}+\lambda \mu_{b,h}-\lambda \mu_{a,h}=0 \quad {\rm a. e.  \ in} \ \mathrm{\Omega},\label{eqn3:KKT for P_lambda,h}\\
       &&(\mu_{a,h},a-\lambda \overline u_{\lambda,h}-\overline y_{\lambda,h})=(\mu_{b,h},\lambda \overline u_{\lambda,h}+\overline y_{\lambda,h}-b)=0,\label{eqn4:KKT for P_lambda,h} \\
       &&\mu_{a,h}(x)\geq0,\ \ \mu_{b,h}(x)\geq0 \quad {\rm a. e.  \ in} \ \mathrm{\Omega}, \label{eqn5:KKT for P_lambda,h} \\
       &&a\leq\lambda\overline u_{\lambda,h}+\overline y_{\lambda,h}\leq b \quad\ {\rm a. e.  \ in} \ \mathrm{\Omega}, \label{eqn6:KKT for P_lambda,h}
\end{eqnarray}
\end{subequations}
where $p_{\lambda,h}$ is the adjoint state and $\mu_{a,h},\ \mu_{b,h}$ are Lagrange multipliers. We consider a sequence of positive real numbers ${\lambda_k}$ tending to zero for $k\rightarrow\infty$. We use $(\mathrm{P_{k}})$ to denote the regularized problems associated to $\lambda_k$ and their solutions are denoted by $(\overline{y}_k,\overline{u}_k)$ with an adjoint state $p_k$ and Lagrange multipliers $\mu_{ak},\ \mu_{bk}$. To begin with, we give the following lemma which focuses on the boundedness of the Lagrangian multipliers. Since upper bound and lower bound exist simultaneously in the problem we consider, the proof of the following lemma encounter some difficulties compared with the situation with only one bound. However, we utilize the fact that at least one of the two multipliers is equal to zero and complete the proof.
\begin{lmm}\label{eqn:miu uniformly bounded}
Under Assumption \ref{assumption slater point}, the sequence of Lagrange multipliers $\{\mu_{bk}\}$ and $\{\mu_{ak}\}$  are uniformly bounded in $L^1(\mathrm{\Omega})$.
\end{lmm}
\begin{proof}
Let $u_1=\min(\widehat{u},0),\ u_2=\max(\widehat{u},0)\in L^2(\mathrm{\Omega})$, then we have $u_1(x)\leq0,\ u_2(x)\geq0 \ \ {\rm a. e. \ in} \ \mathrm{\Omega}$. Then from the maximum principle for the state equation, we have $(Su_1)(x)<b,\ a<(Su_2)(x)\ \forall x\in \overline{\mathrm{\Omega}}$. So $\forall \lambda\geq0$, there exists $\tau_1,\ \tau_2>0$ such that
\begin{equation}\label{Slater point}
\begin{aligned}
&\lambda u_1(x)+(Su_1)(x)\leq b-\tau_1\quad {\rm a. e. \ in} \ \mathrm{\Omega},\\
&a+\tau_2\leq\lambda u_2(x)+(Su_2)(x)\quad {\rm a. e. \ in} \ \mathrm{\Omega}.
\end{aligned}
\end{equation}
Let $\widehat{u}_{1,k}=u_1-\overline{u}_k,\ \widehat{u}_{2,k}=\overline{u}_k-u_2$, then using (\ref{Slater point}) we arrive at
\begin{equation}\label{miu uniformly bounded 1}
\begin{aligned}
&\tau_1+\lambda_k\overline{u}_k(x)+(S\overline{u}_k)(x)-b\leq-(\lambda_k\widehat{u}_{1,k}(x)+(S\widehat{u}_{1,k})(x))\quad {\rm a. e.  \ in} \ \mathrm{\Omega},\\
&\tau_2+a-\lambda_k\overline{u}_k(x)-(S\overline{u}_k)(x)\leq-(\lambda_k\widehat{u}_{2,k}(x)+(S\widehat{u}_{2,k})(x))\quad {\rm a. e.  \ in} \ \mathrm{\Omega}.
\end{aligned}
\end{equation}
We multiply the two formulas in (\ref{miu uniformly bounded 1}) by $\mu_{bk}$ and $\mu_{ak}$ respectively, which implies
\begin{equation}\label{miu uniformly bounded 2}
\begin{aligned}
&\int_{\mathrm{\Omega}}\tau_1 \mu_{bk}dx \leq \int_{\mathrm{\Omega}} -(\lambda_k\widehat{u}_{1,k}+S\widehat{u}_{1,k})\mu_{bk}dx\quad {\rm a. e.  \ in} \ \mathrm{\Omega},\\
&\int_{\mathrm{\Omega}}\tau_2 \mu_{ak}dx \leq \int_{\mathrm{\Omega}} -(\lambda_k\widehat{u}_{2,k}+S\widehat{u}_{2,k})\mu_{ak}dx\quad {\rm a. e.  \ in} \ \mathrm{\Omega}.
\end{aligned}
\end{equation}
Since (\ref{eqn3:KKT for P_lambda}) is equivalent to
\begin{equation}\label{miu uniformly bounded 4}
\int_{\mathrm{\Omega}} (\alpha \overline u_k+G^*(G\overline u_k-y_d+\mu_{bk}-\mu_{ak})+\lambda_k \mu_{bk}-\lambda_k \mu_{ak})z dx=0,\quad \forall z\in L^2(\mathrm{\Omega}).
\end{equation}
We know that at least one of $\mu_{bk}$ and $\mu_{ak}$ is $0$. When $\mu_{ak}$ is $0$, let $z=\widehat{u}_{1,k}$ in (\ref{miu uniformly bounded 4}), then we arrive at
\begin{equation}\label{miu uniformly bounded 5}
\int_{\mathrm{\Omega}} -(\lambda_k\widehat{u}_{1,k}+G\widehat{u}_{1,k})\mu_{bk}dx=\int_{\mathrm{\Omega}} (\alpha\overline u_k+G^*(G\overline u_k-y_d))\widehat{u}_{1,k}dx.
\end{equation}
When $\mu_{bk}$ is $0$, let $z=\widehat{u}_{2,k}$ in (\ref{miu uniformly bounded 4}), then we get
\begin{equation}\label{miu uniformly bounded 6}
\int_{\mathrm{\Omega}} -(\lambda_k\widehat{u}_{2,k}+G\widehat{u}_{2,k})\mu_{ak}dx=\int_{\mathrm{\Omega}} -(\alpha\overline u_k+G^*(G\overline u_k-y_d))\widehat{u}_{2,k}dx.
\end{equation}
Together with (\ref{miu uniformly bounded 2}), we arrive at
\begin{equation}\label{miu uniformly bounded 7}
\begin{aligned}
&\int_{\mathrm{\Omega}}\tau_1 \mu_{bk}dx \leq ((\alpha+\|G\|^2)\|\overline{u}_k\|+\|G\|\|y_d\|)(\|u_1\|+\|\overline{u}_k\|),\\
&\int_{\mathrm{\Omega}}\tau_2 \mu_{ak}dx \leq ((\alpha+\|G\|^2)\|\overline{u}_k\|+\|G\|\|y_d\|)(\|u_2\|+\|\overline{u}_k\|).
\end{aligned}
\end{equation}
From the optimality of $\overline{u}_k$, we know the uniform boundedness of ${\overline{u}_k}$ in $L^2(\mathrm{\Omega})$. So we know that $\{\mu_{bk}\}$ and $\{\mu_{ak}\}$  are uniformly bounded in $L^1(\mathrm{\Omega})$.
\end{proof}

Similarly to Lemma \ref{eqn:miu uniformly bounded}, we can prove the uniform boundedness of $\|\mu_{ak,h}\|_{L^1(\mathrm{\Omega})}$ and $\|\mu_{bk,h}\|_{L^1(\mathrm{\Omega})}$ w.r.t $h,\lambda$ by replacing $S$ by $S_h$ and $G$ by $G_h$.

\begin{thrm}\label{An error estimate uniform in lambda}
Let $(\overline{y}_{\lambda},\overline{u}_{\lambda})$ and $(\overline{y}_{\lambda,h},\overline{u}_{\lambda,h})$ be the solutions of {\rm{(\ref{eqn:modified problem lambda})}} and {\rm{(\ref{eqn:discretized problem lambda})}} respectively, then there exists some $0<h_0\leq1$ such that
\begin{equation*}
\|\overline{u}_\lambda-\overline{u}_{\lambda,h}\|+\|\overline{y}_\lambda-\overline{y}_{\lambda,h}\|_{H^1(\mathrm{\Omega})}\leq Ch^{1-\frac{n}{4}}, \quad \forall\ 0<h\leq h_0
\end{equation*}
holds, where $n$ denotes the dimension of $\mathrm{\Omega}$ and $C>0$ is a positive constant which is independent of $\lambda$.
\end{thrm}
\begin{proof}
Subtracting (\ref{eqn3:KKT for P_lambda}) and (\ref{eqn3:KKT for P_lambda,h}), we get:
\begin{equation}
       \alpha(\overline u_{\lambda}-\overline u_{\lambda,h})+(p_{\lambda}-p_{\lambda,h})+\lambda(\mu_b-\mu_{b,h})-\lambda(\mu_a-\mu_{a,h})=0.
\end{equation}
Multiplying the formula above by $\overline u_{\lambda}-\overline u_{\lambda,h}$, we derive
\begin{equation}\label{eqn:*}
\begin{aligned}
       \alpha\|\overline u_{\lambda}-\overline u_{\lambda,h}\|^2&=(p_{\lambda,h}-p_{\lambda},\overline u_{\lambda}-\overline u_{\lambda,h})-(\lambda(\mu_b-\mu_{b,h}),\overline u_{\lambda}-\overline u_{\lambda,h})\\
       &\quad+(\lambda(\mu_a-\mu_{a,h}),\overline u_{\lambda}-\overline u_{\lambda,h})\\
       &=-(\lambda(\mu_b-\mu_{b,h}),\overline u_{\lambda}-\overline u_{\lambda,h})+(\lambda(\mu_a-\mu_{a,h}),\overline u_{\lambda}-\overline u_{\lambda,h})\\
       &\quad+(p^h-p_{\lambda},\overline u_{\lambda}-\overline u_{\lambda,h})+(p_{\lambda,h}-p^h,\overline u_{\lambda}-\overline u_{\lambda,h}),
\end{aligned}
\end{equation}
where $p^h$ is the solution of
\begin{equation}\label{eqn:p^h}
        (\nabla p^h,\nabla z_h)=(\overline y_{\lambda}-y_d-\mu_a+\mu_b,z_h)\ \ \forall z_h\in Z_h,
\end{equation}
$y^h$ is the solution of
\begin{equation}\label{eqn:y^h}
        (\nabla y^h,\nabla z_h)=(\overline u_{\lambda},z_h)\ \ \forall z_h\in Z_h.
\end{equation}
Let $z_h=y^h-\overline y_{\lambda,h}\in Z_h$ in the formula which we get by subtracting (\ref{eqn2:KKT for P_lambda,h}) and (\ref{eqn:p^h}), then we arrive at
\begin{equation}
       (\nabla (p_{\lambda,h}-p^h),\nabla (y^h-\overline y_{\lambda,h}))=(\overline y_{\lambda,h}-\overline y_{\lambda}+\mu_{b,h}-\mu_b-\mu_{a,h}+\mu_a,y^h-\overline y_{\lambda,h}).
\end{equation}
Similarly, let $z_h=p_{\lambda,h}-p^h\in Z_h$ in the formula which we get by subtracting (\ref{eqn:y^h}) and (\ref{eqn1:KKT for P_lambda,h}), then we derive
\begin{equation}
       (\nabla (y^h-\overline y_{\lambda,h}),\nabla (p_{\lambda,h}-p^h))=(\overline u_\lambda-\overline u_{\lambda,h},p_{\lambda,h}-p^h).
\end{equation}
So we can get
\begin{equation*}
\begin{aligned}
       &(p_{\lambda,h}-p^h,\overline u_\lambda-\overline u_{\lambda,h})=(\overline y_{\lambda,h}-\overline y_{\lambda},y^h-\overline y_{\lambda,h})+\underbrace{(\mu_b,\overline y_{\lambda,h}-y^h)}_I+\underbrace{(\mu_{b,h},y^h-\overline y_{\lambda,h})}_{II}\\
       &\qquad\qquad\qquad\qquad\qquad+\underbrace{(\mu_{a},y^h-\overline y_{\lambda,h})}_{III}+\underbrace{(\mu_{a,h},\overline y_{\lambda,h}-y^h)}_{IV}.
\end{aligned}
\end{equation*}
For the term $I$, since $\overline y_{\lambda,h}\leq b-\lambda \overline u_{\lambda,h}$ and $\mu_b\geq0$, we derive
\begin{equation}\label{eqn:I}
\begin{aligned}
       (\mu_b,\overline y_{\lambda,h}-y^h)&\leq(\mu_b,b-\lambda\overline u_{\lambda,h}-y^h-b+\lambda\overline u_{\lambda}+\overline y_{\lambda})\\
       &=(\mu_b,\lambda(\overline u_{\lambda}-\overline u_{\lambda,h}))+(\mu_b,\overline y_{\lambda}-y^h).
\end{aligned}
\end{equation}
For the term $II$, because of $\overline y_{\lambda}\leq b-\lambda \overline u_{\lambda}$ and $\mu_{b,h}\geq0$, we have
\begin{equation}\label{eqn:II}
\begin{aligned}
       (\mu_{b,h},y^h-\overline y_{\lambda,h})&=(\mu_{b,h},\overline y_{\lambda}-\overline y_{\lambda,h})+(\mu_{b,h},y^h-\overline y_{\lambda})\\
       &\leq(\mu_{b,h},b-\lambda\overline u_{\lambda}-\overline y_{\lambda,h}-b+\lambda\overline u_{\lambda,h}+\overline y_{\lambda,h})+(\mu_{b,h},y^h-\overline y_{\lambda})\\
       &\leq(\mu_{b,h},\lambda(\overline u_{\lambda,h}-\overline u_{\lambda}))+(\mu_{b,h},y^h-\overline y_{\lambda}).
\end{aligned}
\end{equation}
For the term $III$, based on $-\overline y_{\lambda,h}\leq -a+\lambda \overline u_{\lambda,h}$ and $\mu_a\geq0$, we arrive at
\begin{equation}\label{eqn:III}
\begin{aligned}
       (\mu_a,y^h-\overline y_{\lambda,h})&\leq(\mu_a,y^h-a+\lambda\overline u_{\lambda,h}+a-\lambda\overline u_{\lambda}-\overline y_{\lambda})\\
       &=(\mu_a,\lambda(\overline u_{\lambda,h}-\overline u_{\lambda}))+(\mu_a,y^h-\overline y_{\lambda})).
\end{aligned}
\end{equation}
For the term $IV$, following from $-\overline y_{\lambda}\leq -a+\lambda \overline u_{\lambda}$ and $\mu_{a,h}\geq0$, we have
\begin{equation}\label{eqn:IV}
\begin{aligned}
       (\mu_{a,h},\overline y_{\lambda,h}-y^h)&=(\mu_{a,h},\overline y_{\lambda,h}-\overline y_{\lambda})+(\mu_{a,h},\overline y_{\lambda}-y^h)\\
       &\leq(\mu_{a,h},\overline y_{\lambda,h}-a+\lambda\overline u_{\lambda}+a-\lambda\overline u_{\lambda,h}-\overline y_{\lambda,h})+(\mu_{a,h},\overline y_{\lambda}-y^h)\\
       &\leq(\mu_{a,h},\lambda(\overline u_{\lambda}-\overline u_{\lambda,h}))+(\mu_{a,h},\overline y_{\lambda}-y^h).
\end{aligned}
\end{equation}
Inserting (\ref{eqn:I}), (\ref{eqn:II}), (\ref{eqn:III}) and (\ref{eqn:IV}) into (\ref{eqn:*}), we get
\begin{equation*}
\begin{aligned}
       &\alpha\|\overline u_{\lambda}-\overline u_{\lambda,h}\|^2\leq-(\lambda(\mu_b-\mu_{b,h}),\overline u_{\lambda}-\overline u_{\lambda,h})+(\lambda(\mu_a-\mu_{a,h}),\overline u_{\lambda}-\overline u_{\lambda,h})\\
       &\qquad \qquad \qquad\ \ \ +(p^h-p_{\lambda},\overline u_{\lambda}-\overline u_{\lambda,h})+(\overline y_{\lambda,h}-\overline y_\lambda,y^h-\overline y_{\lambda,h})\\
       &\qquad \qquad \qquad\ \ \ +(\lambda(\mu_b-\mu_{b,h}),\overline u_{\lambda}-\overline u_{\lambda,h})-(\lambda(\mu_a-\mu_{a,h}),\overline u_{\lambda}-\overline u_{\lambda,h})\\
       &\qquad \qquad \qquad\ \ \ +(\mu_b,\overline y_{\lambda}-y^h)+(\mu_{b,h},y^h-\overline y_{\lambda})+(\mu_a,y^h-\overline y_{\lambda})+(\mu_{a,h},\overline y_{\lambda}-y^h)\\
       &\qquad \qquad \qquad=-\|\overline y_{\lambda}-\overline y_{\lambda,h}\|^2+(\overline y_{\lambda,h}-\overline y_{\lambda},y^h-\overline y_{\lambda})+(p^h-p_{\lambda},\overline u_{\lambda}-\overline u_{\lambda,h})\\
       &\qquad \qquad \qquad\ \ \ +(\mu_b-\mu_{b,h}-\mu_a+\mu_{a,h},\overline y_{\lambda}-y^h),
\end{aligned}
\end{equation*}
which gives
\begin{equation*}
\begin{aligned}
       &\ \ \ \alpha\|\overline u_{\lambda}-\overline u_{\lambda,h}\|^2+\|\overline y_{\lambda}-\overline y_{\lambda,h}\|^2\\
       &\leq(\overline y_{\lambda,h}-\overline y_{\lambda},y^h-\overline y_{\lambda})+(p^h-p_{\lambda},\overline u_{\lambda}-\overline u_{\lambda,h})+(\mu_b-\mu_{b,h}-\mu_a+\mu_{a,h},\overline y_{\lambda}-y^h)\\
       &\leq\frac{1}{2}\|\overline y_{\lambda}-\overline y_{\lambda,h}\|^2+\frac{1}{2}\|y^h-\overline y_{\lambda}\|^2+\frac{2}{\alpha}\|p^h-p_{\lambda}\|^2+\frac{\alpha}{2}\|\overline u_{\lambda}-\overline u_{\lambda,h}\|^2\\
       &\ \ \ +\|\mu_b-\mu_{b,h}-\mu_a+\mu_{a,h}\|_{L^1(\mathrm{\Omega})}\cdot\|\overline y_{\lambda}-y^h\|_{L^{\infty}(\mathrm{\Omega})}.
\end{aligned}
\end{equation*}
Then we arrive at
\begin{equation*}
\begin{aligned}
       &\ \ \ \alpha\|\overline u_{\lambda}-\overline u_{\lambda,h}\|^2+\|\overline y_{\lambda}-\overline y_{\lambda,h}\|^2\\
       &\leq\|\overline y_{\lambda}-y^h\|^2+\frac{4}{\alpha}\|p_{\lambda}-p^h\|^2+2\|\mu_b-\mu_{b,h}-\mu_a+\mu_{a,h}\|_{L^1(\mathrm{\Omega})}\cdot\|\overline y_{\lambda}-y^h\|_{L^{\infty}(\mathrm{\Omega})}.
\end{aligned}
\end{equation*}
It is shown in \cite{casas19852} that the following formula holds
\begin{equation}
\|p_{\lambda}-p^h\|^2\leq h^{4-n}\left(\|\overline{y}_\lambda-y_d\|^2+\|\mu_a\|^2_{L^1}+\|\mu_b\|^2_{L^1}\right).
\end{equation}
Through standard finite element error estimates and the fact that $\|\overline{u}_\lambda\|$ is bounded independent of $\lambda$ resulting from the optimality of $\overline{u}_\lambda$, we know that $\|\overline y_{\lambda}-y^h\|^2\leq Ch^4$ and $\|\overline y_{\lambda}-y^h\|_{L^{\infty}(\mathrm{\Omega})}\leq Ch^{2-\frac{n}{2}}$. Together with Lemma \ref{eqn:miu uniformly bounded}, we have the following estimation for $\|\overline u_{\lambda}-\overline u_{\lambda,h}\|$,
\begin{equation}
\|\overline u_{\lambda}-\overline u_{\lambda,h}\|^2\leq C(h^4+h^{4-n}+h^{2-\frac{n}{2}}),
\end{equation}
which implies $\|\overline u_{\lambda}-\overline u_{\lambda,h}\|\leq Ch^{1-\frac{n}{4}}$. Then as the proof of Theorem \ref{An error estimate for fixed lambda}, we can get
\begin{equation}
\|\overline y_{\lambda}-\overline y_{\lambda,h}\|_{H^1}\leq C(h+h^{1-\frac{n}{4}}),
\end{equation}
which gives the assertion.

\end{proof}

In addition, if we assume that $\overline u_{\lambda}$ is uniformly bounded in $L^{\infty}(\mathrm{\Omega})$, then from \cite{deckelnick2008numerical} we know that
\begin{equation}
\|\overline y_{\lambda}-y^h\|_{L^{\infty}(\mathrm{\Omega})}\leq Ch^2|\log(h)|^2\|\overline u_{\lambda}\|_{L^{\infty}(\mathrm{\Omega})}.
\end{equation}
Then from the proof of Theorem \ref{An error estimate uniform in lambda} we have
\begin{equation}
\|\overline u_{\lambda}-\overline u_{\lambda,h}\|^2\leq C(h^4+h^{4-n}+h^2|\log(h)|^2),
\end{equation}
which implies the following corollary.
\begin{crllr}\label{error estimates uniform in lambda corollary}
Assume that the sequence of optimal solutions to {\rm{(\ref{eqn:modified problem lambda})}} for $\lambda\downarrow0$, denoted by $\{\overline{u}_\lambda\}$, is uniformly bounded in $L^\infty(\mathrm{\Omega})$, and assume further that the solution of {\rm{(\ref{eqn:state equation})}} satisfies $y\in W^{2,q}(\mathrm{\Omega})$ for all $1\leq q <\infty$ if $u\in L^\infty(\mathrm{\Omega})$. Then the sequence of solutions of {\rm(\ref{eqn:discretized problem lambda})}, denoted by $\{\overline{u}_{\lambda,h}\}$ satisfies
\begin{equation*}
\|\overline{u}_\lambda-\overline{u}_{\lambda,h}\|\leq C\max\{h|\log(h)|,h^{2-\frac{n}{2}}\},\qquad\forall\ 0<h\leq h_0
\end{equation*}
where $n$ denotes the dimension of $\mathrm{\Omega}$ and $C$ is a constant independent of $\lambda$ and $h$.
\end{crllr}

\subsection{Analysis for error estimates}
The main novelty with respect to the error estimates of our paper is that we prove the error order of full discretization is not inferior to that of variational discretization, which has been stated in detail in introduction. The overall error consists of two parts: one arising from the regularization and another caused by the discretization. We know from \cite{prufert2008convergence} that for the error resulted from Lavrentiev-regularization, the following theorem holds
\begin{thrm}\label{regularization error}
Let $(y^*,u^*)$ and $(\overline{y}_{\lambda},\overline{u}_{\lambda})$ be the solutions of {\rm{(\ref{eqn:original problem})}} and {\rm{(\ref{eqn:modified problem lambda})}}, then the following error estimate holds
\begin{equation*}
\|u^*-\overline u_{\lambda}\|\leq c\sqrt{\lambda},
\end{equation*}
where $c$ is a constant independent of $\lambda$.
\end{thrm}

Combining Theorem \ref{regularization error} with Theorem \ref{An error estimate for fixed lambda} and Corollary \ref{error estimates uniform in lambda corollary}, we arrive at the following results for the overall error.
\begin{equation}\label{overall error estimate 1}
\|u^*-\overline u_{\lambda,h}\|\leq C_1\left(\sqrt{\lambda}+\frac{1}{\sqrt{\lambda}}h+\frac{1}{\lambda^2}\left(h^2+\frac{1}{\lambda}h^3+\frac{1}{\lambda^2}h^4\right)\right)
\end{equation}
and
\begin{equation}\label{overall error estimate 2}
\|u^*-\overline u_{\lambda,h}\|\leq C_2\left(\sqrt{\lambda}+\max\{h|\log(h)|,h^{2-\frac{n}{2}}\}\right),
\end{equation}
where $n=2,\ 3$ denotes the dimension of $\mathrm{\Omega}$ and $C_1$, $C_2$ are positive constants independent of $\lambda$ and $h$. As we said in Introduction, the error order of full discretization is not inferior to that of variational discretization. It is clear from (\ref{overall error estimate 1}) and (\ref{overall error estimate 2}) that when $\lambda$ is fixed, both two error estimates decrease as $h$ declines until reaching a lower bound resulting from term $\sqrt{\lambda}$, i.e. Lavrentiev regularization. However, for fixed $h$, the first error estimate may decrease also may increase as $\lambda$ declines because term $\sqrt{\lambda}$ and term $\frac{1}{\sqrt{\lambda}}h+\frac{1}{\lambda^2}\left(h^2+\frac{1}{\lambda}h^3+\frac{1}{\lambda^2}h^4\right)$ exist simultaneously. While the second error estimate may decline until reaching a lower bound also may remain unchanged as $\lambda$ decreases. These statements declare that for fixed $h$, it is not the smaller $\lambda$ the better. Additionally, both (\ref{overall error estimate 1}) and (\ref{overall error estimate 2}) give an upper bound for the error, which one is a better estimate also depends on the values of $C_1$ and $C_2$. So different problems may have various error variation trend. We could verify the statement above through the numerical experiments in Section \ref{sec:5}.

\section{hADMM and two-phase strategy}
\label{sec:4}
The error of utilizing numerical methods to solve PDE constrained problem consists of two parts: discretization error and the error of algorithm for discretized problem. The error order of piecewise linear finite element method is $O(h)$, which makes the discretization error account for the main part. So algorithms of high precision do not make much sense, instead will waste much computations. Thus using heterogeneous ADMM (hADMM), which is a fast and efficient first order algorithm, to get a solution of moderate precision is sufficient. Heterogeneous ADMM is different from the classical ADMM, where two different norms are applied in the first two subproblems. However, in order to satisfy the need for more accurate solution, a two-phase strategy is also presented, in which the PDAS method is used as a postprocessor of the hADMM algorithm. However, we should emphasize that here the `accurate' refers to the KKT precision of the numerical algorithm but not the error between exact solution and numerical solution.

To rewrite the discretized problem into a matrix-vector form, we define the following matrices
\begin{equation}
   {K_h}=\left(\int_{\mathrm{\Omega}_h}\nabla\phi_i\cdot\nabla\phi_j\ dx\right)_{i,j=1}^{N_h}\quad{\rm{and}}\quad {M_h}=\left(\int_{\mathrm{\Omega}_h}\phi_i\cdot\phi_j\ dx\right)_{i,j=1}^{N_h},
\end{equation}
where ${K_h}$ and ${M_h}$ denote the finite element stiffness matrix and mass matrix respectively. Let
\begin{equation}
  y_{d,h}(x)=\sum\limits_{i=1}^{N_h}y_d^i\phi_i(x)
\end{equation}
be the $L^2$-projection of $y_d$ onto $Z_h$, where $y_d^i=y_d(x^i)$. The lump mass matrix $W_h$ is defined by
\begin{equation}\label{Wh}
  W_h={\rm{diag}}\left(\int_{\mathrm{\Omega}_h}\ \phi_i(x)\ dx\right)^{N_h}_{i=1},
\end{equation}
which is a diagonal matrix. Actually, each principal diagonal element of $W_h$ is twice as the counterpart of $M_h$. For the mass matrix $M_h$ and the lump mass matrix $W_h$, the following proposition hold.
\begin{prpstn}{\rm\cite[Table 1]{wathen1987realistic}}\label{eqn:martix properties}
$\forall$ $z\in \mathbb{R}^{N_h}$, the following inequalities hold:
\begin{equation*}
  \|z\|^2_{M_h}\leq\|z\|^2_{W_h}\leq c\|z\|^2_{M_h},\quad where\quad c=\left\{
  \begin{aligned}
  &4\quad if\quad n=2,\\
  &5\quad if\quad n=3.
  \end{aligned}\right.
\end{equation*}
\end{prpstn}

For simplicity, we use the symbol before discretization to denote the column vectors of the coefficients of the functions with respect to the basis $\{\phi_i(x)\} _{i=1}^{N_h}$ which are discretized above, for example, $y=(y_1,y_2,\cdots,y_{N_h})^T\in \mathbb{R}^{N_h}$. Then we can rewrite the problem (\ref{eqn:discretized problem lambda}) and (\ref{eqn:discretized problem vu}) into a matrix-vector form, which are the actual versions we apply the hADMM algorithm and PDAS method to respectively

\begin{equation}\label{eqn:discretized matrix-vector form vu}
\left\{\begin{aligned}
        \min\limits_{y,u,v\in {\mathbb R}^{N_h}}\ &{J}_h(y,u)= \frac{1}{2}\|y-y_d\|_{{M_h}}^{2}+\frac{\alpha}{2}\|u\|_{{M_h}}^{2}\\
        {\rm{s.t.}}\quad\ &K_hy={M_h}u,\\
        &v-\lambda u-y=0, \\
        &v\in[a,b]^{N_h}.
                          \end{aligned} \right.\tag{${\mathrm{\widehat{P}}'}_{\lambda,h}$}
\end{equation}

\begin{equation}\label{eqn:discretized matrix-vector form lambda}
\left\{\begin{aligned}
        \min\limits_{y,u\in {\mathbb R}^{N_h}}\ \ &{J}_h(y,u)= \frac{1}{2}\|y-y_d\|_{{M_h}}^{2}+\frac{\alpha}{2}\|u\|_{{M_h}}^{2}\\
        {\rm{s.t.}}\quad\ &K_hy={M_h}u,\\
        &\lambda u+y\in[a,b]^{N_h}.
                          \end{aligned} \right.\tag{${\mathrm{P}'}_{\lambda , h}$}
\end{equation}

In the process of implementation, if a solution with moderate accuracy is sufficient, hADMM algorithm is applied. In addition, if more accurate solution (`accurate' refers to the KKT precision of the numerical algorithm but not the error between exact solution and numerical solution) is required, a two-phase strategy is employed, in which the PDAS method is used as a postprocessor of the hADMM algorithm. The following two subsections focus on the hADMM algorithm and the PDAS method respectively.

\subsection{Two ADMM-type algorithms for (\ref{eqn:discretized matrix-vector form vu})}
\label{sec:4.1}
Since the stiffness matrix ${K_h}$ and the mass matrix ${M_h}$ are symmetric positive definite matrices, we can rewrite (\ref{eqn:discretized matrix-vector form vu}) into the reduced form
\begin{equation}\label{eqn:reduced discretized matrix-vector form vu}
\left\{\begin{aligned}
        \min\limits_{u,v\in {\mathbb R}^{N_h}}\ \ &{J}_h(y,u)= \frac{1}{2}\|K_h^{-1}{M_h}u-y_d\|_{{M_h}}^{2}+\frac{\alpha}{2}\|u\|_{{M_h}}^{2}\\
        {\rm{s.t.}}\quad \ &v-\lambda u-K_h^{-1}{M_h}u=0, \\
        &v\in[a,b]^{N_h}.
                          \end{aligned} \right.\tag{$\mathrm{R}{\mathrm{\widehat{P}}'}_{\lambda,h}$}
\end{equation}
In order to show the differences between our hADMM and classical ADMM, we give the details of these two algorithms respectively. First, let us focus on classical ADMM.
\subsubsection{Classical ADMM}
We can see from the content below that the first subproblem of classical ADMM has to solve a $3*3$ block equation system. It can be reduced into a $2*2$ block equation system, however, it will introduce additional computation of $M_h^{-1}$. More importantly, classical ADMM algorithm is not mesh independent.

The augmented Lagrangian function of (\ref{eqn:reduced discretized matrix-vector form vu}) is:
\begin{equation}
\begin{aligned}
        L_{\sigma}(v,u;\mu)=&\frac{1}{2}\|K_h^{-1}{M_h}u-y_d\|_{{M_h}}^{2}+\frac{\alpha}{2}\|u\|_{{M_h}}^{2}+(\mu,v-\lambda u-K_h^{-1}{M_h}u)\\
        &+\frac{\sigma}{2}\|v-\lambda u-K_h^{-1}{M_h}u\|^2+\delta_{[a,b]^{N_h}}(v),
\end{aligned}
\end{equation}
where $\mu \in {\mathbb R}^{N_h}$ is the Lagrange multiplier and $\sigma>0$ is a penalty parameter. We give the three main steps at $k$-th iteration.
\begin{equation*}\left\{
\begin{aligned}
        step 1: u^{k+1}&=\arg\min_u L_{\sigma}(v^k,u;\mu^k)\\
        step 2: v^{k+1}&=\arg\min_v L_{\sigma}(v,u^{k+1};\mu^k)\\
        step 3: \mu^{k+1}&= \mu^k+\sigma(v^{k+1}-\lambda u^{k+1}-y^{k+1})
\end{aligned}\right.
\end{equation*}
Now let us give the details about two subproblems with respect to $u$ and $v$ respectively. The first subproblem is equivalent to the following problem
\begin{equation}\label{cADMM_step1}
\begin{aligned}
\min_{y,u\in{\mathbb R}^{N_h}}\quad  &\frac{1}{2}\|y-y_d\|_{{M_h}}^{2}+\frac{\alpha}{2}\|u\|_{{M_h}}^{2}+(\mu^k,v^k-\lambda u-y)+\frac{\sigma}{2}\|v^k-\lambda u-y\|^2\\
{\rm s.t.}\qquad &K_hy-M_hu=0,
\end{aligned}
\end{equation}
whose Lagrangian function is
\begin{equation*}
L_1(y,u;p)=\frac{1}{2}\|y-y_d\|_{{M_h}}^{2}+\frac{\alpha}{2}\|u\|_{{M_h}}^{2}+(\mu^k,v^k-\lambda u-y)+\frac{\sigma}{2}\|v^k-\lambda u-y\|^2+(p,K_hy-M_hu),
\end{equation*}
where $p$ is the Lagrangian multiplier corresponding to the equality constraint $K_hy-M_hu=0$. Then the KKT conditions of (\ref{cADMM_step1}) are
\begin{equation*}\left\{
\begin{aligned}
&M_h(y-y_d)-\mu^k-\sigma(v^k-\lambda u-y)+{K_h}^Tp=0\\
&\alpha M_h u-\lambda\mu^k-\lambda\sigma(v^k-\lambda u-y)-{M_h}^Tp=0\\
&K_hy-M_hu=0
\end{aligned}\right.
\end{equation*}
\begin{equation}\label{cADMM step1 es}
\Leftrightarrow\left[
  \begin{array}{ccc}
    {M_h}+\sigma I & \lambda\sigma I & {K_h}^T \\
    \lambda\sigma I & \lambda^2\sigma I+\alpha{M_h} & -{M_h}^T \\
    K_h & -{M_h} & 0 \\
  \end{array}
\right]\left[
         \begin{array}{c}
           y^{k+1} \\
           u^{k+1} \\
           p^{k+1} \\
         \end{array}
       \right]=\left[
                 \begin{array}{c}
                   {M_h}y_d+\mu^k+\sigma v^k \\
                   \lambda(\mu^k+\sigma v^k) \\
                   0 \\
                 \end{array}
               \right].
\end{equation}

The second subproblem is equivalent to the following problem
\begin{equation}\label{cADMM_step2}
\begin{aligned}
\min_{v\in{\mathbb R}^{N_h}}\quad  &(\mu,v-\lambda u^{k+1}-y^{k+1})+\frac{\sigma}{2}\|v-\lambda u^{k+1}-y^{k+1}\|^2\\
{\rm s.t.}\quad\ &v\in[a,b]^{N_h},
\end{aligned}
\end{equation}
whose object function is a quadratic function, so it has a closed form solution
\begin{equation}
v^{k+1}={\rm\Pi}_{[a,b]^{N_h}}\left(\lambda u^{k+1}+y^{k+1}-\frac{\mu^k}{\sigma}\right).
\end{equation}


\subsubsection{Heterogeneous ADMM (hADMM)}
The essential difference between hADMM and classical ADMM is that the former adopts two different weighted norms in two subproblems in each iteration. It is clear from the content below that the first subproblem of hADMM only has to solve a $2*2$ block system without any additional computations, which can be solved by generalized minimal residual (GMRES) with preconditioning matrix, and the second subproblem has a closed form solution. More importantly, The numerical results in Section \ref{sec:5} indicate that our hADMM algorithm is mesh independent, while classical ADMM is not. Additionally, to construct the relation between the continuous problem and discretized problem, proposing hADMM algorithm is a natural idea. Following the hADMM proposed  in \cite{song2016two}, whose idea is to employ two different weighted norms in two subproblems, the weighted augmented Lagrangian function of (\ref{eqn:reduced discretized matrix-vector form vu}) is:
\begin{equation}
\begin{aligned}
        \widetilde{L}_{\sigma}(v,u;\mu)=&\frac{1}{2}\|K_h^{-1}{M_h}u-y_d\|_{{M_h}}^{2}+\frac{\alpha}{2}\|u\|_{{M_h}}^{2}+(\mu,v-\lambda u-K_h^{-1}{M_h}u)_{M_h}\\
        &+\frac{\sigma}{2}\|v-\lambda u-K_h^{-1}{M_h}u\|_{M_h}^2+\delta_{[a,b]^{N_h}}(v),
\end{aligned}
\end{equation}
where $\mu \in {\mathbb R}^{N_h}$ is the Lagrange multiplier and $\sigma>0$ is a penalty parameter. The three steps in each iteration of hADMM algorithm are as follows
\begin{equation*}\left\{
\begin{aligned}
        step 1: u^{k+1}&=\arg\min_u \widetilde{L}_{\sigma}(v^k,u;\mu^k)\\
        step 2: v^{k+1}&=\arg\min_v \widetilde{L}_{\sigma}(v,u^{k+1};\mu^k)\\
        step 3: \mu^{k+1}&= \mu^k+\sigma(v^{k+1}-\lambda u^{k+1}-y^{k+1})
\end{aligned}\right.
\end{equation*}
Now let us give the details about two subproblems with respect to $u$ and $v$ respectively. The first subproblem is equivalent to the following problem
\begin{equation}\label{hADMM_step1}
\begin{aligned}
\min_{y,u\in{\mathbb R}^{N_h}}\quad  &\frac{1}{2}\|y-y_d\|_{{M_h}}^{2}+\frac{\alpha}{2}\|u\|_{{M_h}}^{2}+(\mu^k,v^k-\lambda u-y)_{M_h}+\frac{\sigma}{2}\|v^k-\lambda u-y\|_{M_h}^2\\
{\rm s.t.}\qquad &K_hy-M_hu=0,
\end{aligned}
\end{equation}
whose Lagrangian function is
\begin{equation*}
L_2(y,u;p)=\frac{1}{2}\|y-y_d\|_{{M_h}}^{2}+\frac{\alpha}{2}\|u\|_{{M_h}}^{2}+(\mu^k,v^k-\lambda u-y)_{M_h}+\frac{\sigma}{2}\|v^k-\lambda u-y\|_{M_h}^2+(p,K_hy-M_hu),
\end{equation*}
where $p$ is the Lagrangian multiplier corresponding to the equality constraint $K_hy-M_hu=0$. Since the smoothness of (\ref{hADMM_step1}), solving it is equivalent to solving the following linear system
\begin{equation}\label{hADMM step1 es}
\left[
  \begin{array}{ccc}
    (1+\sigma) {M_h} & \lambda\sigma {M_h} & K_h^T \\
    \lambda\sigma {M_h} & (\lambda^2\sigma+\alpha){M_h} & -{M_h}^T \\
    K_h & -{M_h} & 0 \\
  \end{array}
\right]\left[
         \begin{array}{c}
           y^{k+1} \\
           u^{k+1} \\
           p^{k+1} \\
         \end{array}
       \right]=\left[
                 \begin{array}{c}
                   {M_h}(y_d+\mu^k+\sigma v^k) \\
                   \lambda {M_h}(\mu^k+\sigma v^k) \\
                   0 \\
                 \end{array}
               \right],
\end{equation}
from which we derive that
\begin{equation}
u^{k+1}=\frac{1}{\lambda^2\sigma+\alpha}(p^{k+1}-\lambda\sigma y^{k+1}+\lambda(\mu^{k}+\sigma v^{k})).
\end{equation}
Then (\ref{hADMM step1 es}) could be reduced into the following equation system without any additional calculation.
\begin{equation}\label{as2*2}
\left[
  \begin{array}{cc}
    (1+\frac{\sigma\alpha}{\lambda^2\sigma+\alpha}) {M_h} & \frac{\lambda\sigma}{\lambda^2\sigma+\alpha} {M_h}+K_h^T \\
    -\frac{\lambda\sigma}{\lambda^2\sigma+\alpha} {M_h}-K_h & \frac{1}{\lambda^2\sigma+\alpha}{M_h} \\
  \end{array}
\right]\left[
         \begin{array}{c}
           y^{k+1} \\
           p^{k+1} \\
         \end{array}
       \right]
       =\left[
                 \begin{array}{c}
                   {M_h}y_d+\frac{\alpha}{\lambda^2\sigma+\alpha}{M_h}(\mu^k+\sigma v^k) \\
                   -\frac{\lambda}{\lambda^2\sigma+\alpha}{M_h}(\mu^k+\sigma v^k) \\
                 \end{array}
               \right].
\end{equation}
It is seen that the hADMM only has to solve a $2*2$ block equation system in the first subproblem in each iteration. We should emphasize here that writing the optimality conditions in the form of an antisymmetric matrix can make it more convenient for the design of the preconditioning matrix and the equation system can be solved by GMRES with preconditioner. (\ref{as2*2}) can also be written into a symmetric matrix, however, some of the principle elements of the coefficient matrix will be negative. Utilizing PCG or MINRES to silve it will not have advantages than solving (\ref{as2*2}) by GMRES.

The second subproblem is equivalent to the following problem
\begin{equation}\label{cADMM_step2}
\begin{aligned}
\min_{v\in{\mathbb R}^{N_h}}\quad  &(\mu,v-\lambda u^{k+1}-y^{k+1})_{M_h}+\frac{\sigma}{2}\|v-\lambda u^{k+1}-y^{k+1}\|_{M_h}^2\\
{\rm s.t.}\quad\ &v\in[a,b]^{N_h},
\end{aligned}
\end{equation}
which does not have a closed form solution, we replace the term $\frac{\sigma}{2}\|v-\lambda u^{k+1}-y^{k+1}\|^2_{M_h}$ by $\frac{\sigma}{2}\|v-\lambda u^{k+1}-y^{k+1}\|^2_{W_h}$, where $W_h$ is the lump mass matrix defined in (\ref{Wh}). Then the second subproblem is transformed to the following optimization problem
\begin{equation}
\begin{aligned}
\min_{v\in{\mathbb R}^{N_h}}\quad &(\mu^k,v-\lambda u^{k+1}-y^{k+1})_{M_h}+\frac{\sigma}{2}\|v-\lambda u^{k+1}-y^{k+1}\|_{W_h}^2\\
{\rm s.t.}\quad\ &v\in[a,b]^{N_h},
\end{aligned}
\end{equation}
whose solution has the following closed form
\begin{equation}
v^{k+1}={\rm\Pi}_{[a,b]^{N_h}}\left(\lambda u^{k+1}+y^{k+1}-\frac{{W_h}^{-1}{M_h}\mu^k}{\sigma}\right).
\end{equation}
Although this will introduce the computation of $W_h^{-1}$, $W_h$ is a diagonal matrix, whose inverse will not cost much computation.

Based on the content above, we give the frame of the hADMM algorithm:
\begin{algorithm}[H]\footnotesize
\caption{heterogeneous ADMM (hADMM) algorithm for (\ref{eqn:reduced discretized matrix-vector form vu})}
\label{algo:hADMM for problem RPVU_h}
\leftline{\qquad Initialization: Give initial point $(v^0, \mu^0)\in \mathbb{R}^{N_h} \times \mathbb{R}^{N_h} $ and a tolerant parameter $\tau >0$. Set $k=0$.}
\begin{description}
\item[\textbf{Step 1}] Compute $(y^{k+1},u^{k+1})$ through solving the following equation system
\begin{equation*}
\left[
  \begin{array}{cc}
    (1+\frac{\sigma\alpha}{\lambda^2\sigma+\alpha}) {M_h} & \frac{\lambda\sigma}{\lambda^2\sigma+\alpha} {M_h}+K_h^T \\
    -\frac{\lambda\sigma}{\lambda^2\sigma+\alpha} {M_h}-K_h & \frac{1}{\lambda^2\sigma+\alpha}{M_h} \\
  \end{array}
\right]\left[
         \begin{array}{c}
           y^{k+1} \\
           p^{k+1} \\
         \end{array}
       \right]
       =\left[
                 \begin{array}{c}
                   {M_h}y_d+\frac{\alpha}{\lambda^2\sigma+\alpha}{M_h}(\mu^k+\sigma v^k) \\
                   -\frac{\lambda}{\lambda^2\sigma+\alpha}{M_h}(\mu^k+\sigma v^k) \\
                 \end{array}
               \right].
\end{equation*}
Compute $u^{k+1}$ as follows
\begin{equation*}
u^{k+1}=\frac{1}{\lambda^2\sigma+\alpha}(p^{k+1}-\lambda\sigma y^{k+1}+\lambda(\mu^{k}+\sigma v^{k})).
\end{equation*}
\item[\textbf{Step 2}] Compute $v^{k+1}$ as follows
       \begin{eqnarray*}
       v^{k+1}&=&{\rm\Pi}_{[a,b]^{N_h}}\left(\lambda u^{k+1}+y^{k+1}-\frac{{W_h}^{-1}{M_h}\mu^k}{\sigma}\right).
       \end{eqnarray*}
\item[\textbf{Step 3}] Compute $\mu^{k+1}$ as follows
\begin{eqnarray*}
    \mu^{k+1} &=& \mu^k+\sigma(v^{k+1}-\lambda u^{k+1}-y^{k+1}).
\end{eqnarray*}
\item[\textbf{Step 4}] If a termination criterion is met, Stop; else, set $k:=k+1$ and go to Step 1.
\end{description}
\end{algorithm}

For the convergence result of the heterogeneous ADMM algorithm, we have the following theorem.
\begin{thrm}{\rm\cite[Theorem 4.5]{song2016two}}
Let $(y^*,u^*,v^*,p^*,\mu^*)$ be the KKT point of {\rm(\ref{eqn:discretized matrix-vector form vu})}. $\{(u^k,v^k,\mu^k)\}$ is generated by Algorithm \ref{algo:hADMM for problem RPVU_h} with the associated state $\{y^k\}$ and adjoint state $\{p^k\}$, then we have 
\begin{eqnarray*}
  &\lim\limits_{k\rightarrow\infty}^{}\{\|u^k-u^*\|+\|v^k-v^*\|+\|\mu^k-\mu^*\|\}=0,\\
  &\lim\limits_{k\rightarrow\infty}^{}\{\|y^k-y^*\|+\|p^k-p^*\|\}=0.
\end{eqnarray*}
\end{thrm}

\subsection{Primal-Dual Active Set method as postprocessor}
\label{sec:4.2}
As we have said above, the error of utilizing numerical methods to solve PDE constrained problem consists of two parts: discretization error and the error of algorithm for discretized problem, in which the discretization error account for the main part. Algorithms of high precision do not make much sense but waste computations in practice. In general, using hADMM algorithm to get a solution of moderate precision is sufficient. Although algorithms of high precision are not necessary, we also provide a two-phase strategy to satisfy the requirement for numerical solution of high precision, in which the PDAS method is used as a postprocessor of the hADMM algorithm. The PDAS method was used to solve control constrained elliptic optimal control problem in \cite{bergounioux2002primal}. In \cite{hintermuller2002primal}, the authors show its relation to semismooth Newton method, which can be used to prove its local superlinear convergence. We employ the PDAS method to (\ref{eqn:discretized matrix-vector form lambda}), whose Lagrangian function is:
\begin{equation*}
        \widehat{L}(v,u;\mu)=\frac{1}{2}\|y-y_d\|_{{M_h}}^{2}+\frac{\alpha}{2}\|u\|_{{M_h}}^{2}+(p,K_hy-M_hu)+(\mu_a,a-\lambda u-y)+(\mu_b,\lambda u+y-b),
\end{equation*}
where $\mu_a,\mu_b \in {\mathbb R}^{N_h}$ are the Lagrange multipliers. Then the KKT conditions of (\ref{eqn:discretized matrix-vector form lambda}) are
\begin{equation}\label{KKT mu=mu_b-mu_a}
           \left\{ \begin{aligned}
        &{M_h}(y-y_d)+K_h^Tp-\mu_a+\mu_b=0,\\
        &\alpha {M_h}u-{M_h}^Tp-\lambda\mu_a+\lambda\mu_b=0,\\
        &K_hy-{M_h}u=0,\\
        &\mu_a\geq0,\quad a-\lambda u-y\leq0,\quad (\mu_a,a-\lambda u-y)=0,\\
        &\mu_b\geq0,\quad \lambda u+y-b\leq0,\quad (\mu_b,\lambda u+y-b)=0,
        \end{aligned} \right.
\end{equation}
which can be equivalently rewritten as
\begin{equation}\label{KKT mu=mu_b-mu_a}
           \left\{ \begin{aligned}
        &{M_h}(y-y_d)+K_h^Tp-\mu_a+\mu_b=0,\\
        &\alpha {M_h}u-{M_h}^Tp-\lambda\mu_a+\lambda\mu_b=0,\\
        &K_hy-{M_h}u=0,\\
        &\min(\mu_a,\lambda u+y-a)=\mu_a+\min(0,\lambda u+y-a-\mu_a)=0\\
        &\min(\mu_b,b-\lambda u-y)=\mu_b+\min(0,b-\lambda u-y-\mu_b)=0.
        \end{aligned} \right.
\end{equation}
Let $\mu=\mu_b-\mu_a$, then (\ref{KKT mu=mu_b-mu_a}) can be reduced into the following $4*4$ block system
\begin{equation}\label{KKT mu=mu_b-mu_a}
           \left\{ \begin{aligned}
        &{M_h}(y-y_d)+K_h^Tp+\mu=0,\\
        &\alpha {M_h}u-{M_h}^Tp+\lambda\mu=0,\\
        &K_hy-{M_h}u=0,\\
        &\mu-\max(0,\mu+\lambda u+y-b)-\min(0,\lambda u+y-a-\mu)=0.
        \end{aligned} \right.
\end{equation}
We define the active and inactive sets as
\begin{eqnarray}
  \mathcal{A}_{a,h} &=& \{i\in \{1,2,...,N_h\}: \lambda u_i+y_i+\mu_{i}-a<0\}, \\
  \mathcal{A}_{b,h} &=& \{i\in \{1,2,...,N_h\}: \lambda u_i+y_i+\mu_{i}-b>0\}, \\
  \mathcal{I}&=& \{1,2,\cdots,N_h\}\backslash(\mathcal{A}_{a,h} \cup \mathcal{A}_{b,h})
\end{eqnarray}
and note that the following properties hold
\begin{equation}
\begin{aligned}
  &\lambda u_i+y_i=a\quad \mathrm{on}\ \mathcal{A}_{a,h},\qquad \lambda u_i+y_i=b\quad \mathrm{on}\ \mathcal{A}_{b,h},\\
  &\mu_{i}<0\ \ \mathrm{on}\ \mathcal{A}_{a,h},\ \ \mu_{i}>0\quad \mathrm{on}\ \mathcal{A}_{b,h},\ \ \mu_{i}=0\quad \mathrm{on}\ \mathcal{I}.
\end{aligned}
\end{equation}
Let \begin{equation}\label{EaEb}
           (E_a)_{ij}=\left\{ \begin{aligned}
        &1\qquad i=j\ {\rm{and}}\ i\in \mathcal{A}_{a,h},\\
        &0\qquad {\rm{else}},
        \end{aligned} \right.
        \quad (E_b)_{ij}=\left\{ \begin{aligned}
        &1\qquad i=j\ {\rm{and}}\ i\in \mathcal{A}_{b,h},\\
        &0\qquad {\rm{else}},
        \end{aligned} \right.
\end{equation}
then we can rewrite the optimal system (\ref{KKT mu=mu_b-mu_a}) into a linear system
\begin{equation}
\left[
  \begin{array}{cccc}
    {M_h} & 0 & E_a+E_b & K_h^T \\
    0 & \alpha {M_h} & \lambda (E_a+E_b) & -{M_h}^T \\
    E_a+E_b & \lambda (E_a+E_b) & I-E_a-E_b & 0 \\
    {K_h} & -{M_h} & 0 & 0 \\
  \end{array}
\right]\left[
         \begin{array}{c}
           y \\
           u \\
           \mu \\
           p \\
         \end{array}
       \right]=\left[
                 \begin{array}{c}
                   {M_h}y_d \\
                   0 \\
                   E_aa+E_bb \\
                   0 \\
                 \end{array}
               \right].
\end{equation}
It is shown in \cite{bergounioux2002primal} that whether the two consecutive active sets equal is a termination criterion for the primal-dual active set method. Following the content above, we give the frame of the PDAS method:
\begin{algorithm}[H]\footnotesize
\caption{\ Primal-Dual Active Set (PDAS) algorithm for (\ref{eqn:discretized matrix-vector form lambda})}
\label{Primal-Dual Active Set (PDAS)}
\leftline{\qquad Initialization: Choose initial point $y^0$, $u^0$, $p^0$ and $\mu^0\in{\mathbb R}^{N_h}$; Set $k=0$.}
\begin{description}
\item[\textbf{Step 1}] Determine the following subsets of $\{1,2,...,N_h\}$ (Active and Inactive sets)
\begin{eqnarray*}
\mathcal{A}^{k}_{a,h} &=& \{i\in \{1,2,...,N_h\}: \lambda u^{k}_i+y^{k}_i+\mu^{k}_{i}-a<0\}, \\
\mathcal{A}^{k}_{b,h}&=& \{i\in \{1,2,...,N_h\}: \lambda u^{k}_i+y^{k}_i+\mu^{k}_{i}-b>0\}, \\
\mathcal{I}^{k}&=& \{1,2,...,N_h\}\backslash (\mathcal{A}^{k}_{a,h} \cup \mathcal{A}^{k}_{b,h}).
\end{eqnarray*}
\item[\textbf{Step 2}] Determine $E_a^{k}$ and $E_b^{k}$ through (\ref{EaEb}) and solve the following system
\begin{equation*}
\left[
  \begin{array}{cccc}
    {M_h} & 0 & E_a^{k}+E_b^{k} & K_h^T \\
    0 & \alpha {M_h} & \lambda (E_a^{k}+E_b^{k}) & -{M_h}^T \\
    E_a^{k}+E_b^{k} & \lambda (E_a^{k}+E_b^{k}) & I-E_a^{k}-E_b^{k} & 0 \\
    {K_h} & -{M_h} & 0 & 0 \\
  \end{array}
\right]\left[
         \begin{array}{c}
           y^{k+1} \\
           u^{k+1} \\
           \mu^{k+1} \\
           p^{k+1} \\
         \end{array}
       \right]
       =\left[
                 \begin{array}{c}
                   {M_h}y_d \\
                   0 \\
                   E_a^{k}a+E_b^{k}b \\
                   0 \\
                 \end{array}
               \right].
\end{equation*}
\item[\textbf{Step 3}] If $k>1$, $\mathcal{A}^{k+1}_{a,h}=\mathcal{A}^{k}_{a,h}$ and $\mathcal{A}^{k+1}_{b,h}=\mathcal{A}^{k}_{b,h}$ or a termination criterion is met, Stop; else, set $k:=k+1$ and go to Step 1.
\end{description}
\end{algorithm}

For the convergence result of the PDAS method, we have the following theorem. For more details, we refer to \cite{ulbrich2002nonsmooth, ulbrich2002semismooth, eckstein1992douglas}.
\begin{thrm}
Let ${(u^k,y^k)}$ be generated by Algorithm 2, if the initialization $(u^0,y^0)$ is sufficiently close to the solution $(u^*,y^*)$ of {\rm(\ref{eqn:discretized matrix-vector form lambda})}, then ${(u^k,y^k)}$ converge superlinearly to $(u^*,y^*)$.
\end{thrm}

\section{Numerical Result}
\label{sec:5}
In this section, two numerical experiments are considered. All calculations were performed using MATLAB (R2013a) on a PC with Intel (R) Core (TM) i7-4790K CPU (4.00GHz), whose operation system is 64-bit Windows 7.0 and RAM is 16.0 GB.

In the hADMM algorithm, the accuracy of a numerical solution is measured by the following residual
\begin{equation}
\eta_{\mathrm{A}}=\max\{r_1,\ r_2,\ r_3,\ r_4,\ r_5\},
\end{equation}
where
\begin{equation*}
\begin{aligned}
&r_1=\|M_h(y-y_d)+K_hp-M_h\mu\|,\\
&r_2=\|\alpha M_hu-M_hp-\lambda M_h\mu\|,\\
&r_3=\|v-{\rm\Pi}_{[a,b]}(v-M_h\mu)\|,\\
&r_4=\|K_hy-M_hu\|,\\
&r_5=\|v-\lambda u-y\|_{M_h}.
\end{aligned}
\end{equation*}
Similarly, in the PDAS method, the accuracy of a numerical solution is measured by
\begin{equation}
\eta_{\mathrm{P}}=\max\{\gamma_1,\ \gamma_2,\ \gamma_3,\ \gamma_4\},
\end{equation}
where
\begin{equation*}
\begin{aligned}
&\gamma_1=\|M_h(y-y_d)+K_hp+\mu_{ab}\|,\\
&\gamma_2=\|\alpha M_hu-M_hp+\lambda \mu_{ab}\|,\\
&\gamma_3=\|K_hy-M_hu\|,\\
&\gamma_4=\|\mu_{ab}-\max\{0,\mu_{ab}+\lambda u+y-b\}-\min\{0,\lambda u+y-a+\mu_{ab}\}\|.
\end{aligned}
\end{equation*}
Let $\epsilon$ be a given accuracy tolerance, then the terminal condition is $\eta_{\mathrm{A}}(\eta_{\mathrm{P}})<\epsilon$.

In both two examples, hADMM algorithm and two-phase strategy are employed to get numerical solutions of different precision respectively, i.e. the iteration is terminated with different $\epsilon$. Their convergence behavior are both compared with PDAS method, which is a special semi-smooth Newton method (see \cite{hintermuller2002primal}). There are three tables in both two examples. The first one in each example gives the $L^2$ error of the control, while the last two tables focus on the convergence behavior, including the times of iteration, residual $\eta$ and time, of the hADMM algorithm and the two-phase strategy compared with the PDAS method respectively. In the last two tables, `$\#$dofs' denotes the dimension of the control variable on each grid level, `iter' represents the times of iteration and `residual' represents the precision $\eta$ of the numerical algorithm, which is defined above. In Table \ref{table3} and Table \ref{table6}, two sub columns in the column of `two-phase strategy' record the convergence behavior of two phases, i.e. hADMM and PDAS, respectively.

\begin{xmpl}\label{example:1}
We consider $\mathrm{\Omega}=B_{\frac{5}{2}}(0)$ as the test domain and set $a=-1,\ b=1,\ \alpha=10^{-3}$ and $\sigma=11$ in the first example. The desired state is defined by
\begin{equation*}
y_d(r)=\left\{ \begin{aligned}
        &\ 2\qquad 0\leq r\leq1,\\
        &-2\quad 1< r\leq2,\\
        &\ 0\qquad 2< r\leq2.5.
        \end{aligned} \right.
\end{equation*}

When the exact solution is not known, using numerical solution as relative exact solution is a common method. For more details, one can see \cite{Hinze2009Optimization}. In our practice implementation, we choose $h=\frac{2.5\sqrt{2}}{2^8}$ and $\lambda=10^{-6}$. When $h=\frac{2.5\sqrt{2}}{2^8}$, the scale of data is 306305, which results in a large scale discretized problem. When lambda is too small, the problem will be ill-conditioned and the error will increase on the contrary from the error analysis in Section 3. Through testing with different lambda, e.g. $\lambda=10^{-5.5}$, $10^{-6.5}$ and $10^{-7}$, we find that $\lambda=10^{-6}$ is an appropriate choice. We give the $L^2$ errors $\|u_r^*-\overline{u}_{\lambda,h}\|$ on grids of different sizes with nine different values of $\lambda$ from $10^{-2}$ to $10^{-6}$ in Table \ref{table1}. As an example, the figures of the desired state $y_d$, the numerical state $y_{\lambda,h}$ and numerical control $u_{\lambda,h}$ on the grid of size $h=\frac{2.5\sqrt{2}}{2^6}$ with $\lambda=10^{-4.5}$ are displayed in Figure \ref{example1:yd} and Figure \ref{example1:numerical y and u}. If a solution with moderate accuracy is enough, hADMM algorithm is employed and compared with PDAS method. Both two algorithms are terminated when $\eta_{\mathrm{A}}(\eta_{\mathrm{P}})<10^{-2}$ in this case and the corresponding numerical results are displayed in Table \ref{table2}. In addition, if more accurate solution is required, we employ the two-phase strategy and compare it with PDAS method. In this case, both two algorithms are terminated when $\eta_{\mathrm{A}}(\eta_{\mathrm{P}})<10^{-13}$ and the numerical results are shown in Table \ref{table3}.

From Table \ref{table1}, we can see that for fixed $\lambda$, the error decreases as $h$ declines at first until it reaches a bound resulted from regularization. When $h$ is fixed, the error declines as $\lambda$ decreases generally, while the error shows a rising trend with the last few values of $\lambda$. The numerical results in Table \ref{table1} declares that for fixed $h$ error may increase as $\lambda$ decreases, which verify the error estimates in Section \ref{sec:3}. Table \ref{table2} and Table \ref{table3} show that both the hADMM algorithm and the two-phase strategy are much faster than PDAS method especially when the finite element grid size $h$ is very small. The numerical results in the last two tables verify the efficiency of the hADMM algorithm and the two-phase strategy. We think that the efficiency of hADMM and two-phase strategy will be more obviously when the finite element grid size $h$ get smaller. In our numerical experiment, we think that the finite discretization is fine enough since the dimension of variables of the finest grid level has reached $306305$ and $261121$ in two examples respectively.

\begin{table}[H]\footnotesize
\caption{The $L^2$ error $\|u_r^*-\overline{u}_{\lambda,h}\|$ for Example \ref{example:1}.}
\label{table1}
\begin{center}
\begin{tabular}{@{\extracolsep{\fill}}ccccccccccc}
\hline
\multirow{3}{*}{$\qquad\lambda$} &\multirow{3}{*}{$10^{-2}$} &\multirow{3}{*}{$10^{-2.5}$} &\multirow{3}{*}{$10^{-3}$} &\multirow{3}{*}{$10^{-3.5}$}
&\multirow{3}{*}{$10^{-4}$} &\multirow{3}{*}{$10^{-4.5}$} &\multirow{3}{*}{$10^{-5}$} &\multirow{3}{*}{$10^{-5.5}$} &\multirow{3}{*}{$10^{-6}$}\\
&&&&&&&&&\\
$h$&&&&&&&&&\\
\hline
\multirow{3}{*}{$\frac{2.5\sqrt{2}}{2^4}$} &\multirow{3}{*}{9.8613}&\multirow{3}{*}{5.2420}&\multirow{3}{*}{4.7054}&\multirow{3}{*}{4.7211}&\multirow{3}{*}{4.7837}&\multirow{3}{*}{4.8108}&\multirow{3}{*}{4.8215}&\multirow{3}{*}{4.8250}&\multirow{3}{*}{4.8261}\\
&&&&&&&&&\\
\multirow{3}{*}{$\frac{2.5\sqrt{2}}{2^5}$} &\multirow{3}{*}{9.6385}&\multirow{3}{*}{4.2513}&\multirow{3}{*}{1.9925}&\multirow{3}{*}{1.5824}&\multirow{3}{*}{1.6055}&\multirow{3}{*}{1.6280}&\multirow{3}{*}{1.6360}&\multirow{3}{*}{1.6385}&\multirow{3}{*}{1.6391}\\
&&&&&&&&&\\
\multirow{3}{*}{$\frac{2.5\sqrt{2}}{2^6}$} &\multirow{3}{*}{9.6298}&\multirow{3}{*}{4.2403}&\multirow{3}{*}{1.8113}&\multirow{3}{*}{0.7893}&\multirow{3}{*}{0.5401}&\multirow{3}{*}{0.5360}&\multirow{3}{*}{0.5438}&\multirow{3}{*}{0.5471}&\multirow{3}{*}{0.5484}\\
&&&&&&&&&\\
\multirow{3}{*}{$\frac{2.5\sqrt{2}}{2^7}$} &\multirow{3}{*}{9.5821}&\multirow{3}{*}{4.2302}&\multirow{3}{*}{1.7991}&\multirow{3}{*}{0.7626}&\multirow{3}{*}{0.3373}&\multirow{3}{*}{0.2073}&\multirow{3}{*}{0.1956}&\multirow{3}{*}{0.1997}&\multirow{3}{*}{0.2019}\\
&&&&&&&&&\\
\multirow{3}{*}{$\frac{2.5\sqrt{2}}{2^8}$} &\multirow{3}{*}{9.5361}&\multirow{3}{*}{4.1992}&\multirow{3}{*}{1.7624}&\multirow{3}{*}{0.7588}&\multirow{3}{*}{0.3188}&\multirow{3}{*}{0.1367}&\multirow{3}{*}{0.0660}&\multirow{3}{*}{0.0229}&\multirow{3}{*}{-}\\
&&&&&&&&&\\
&&&&&&&&&\\
\hline
\end{tabular}
\end{center}
\end{table}

\begin{figure}[H]
\centering
\includegraphics[width=0.4\textwidth]{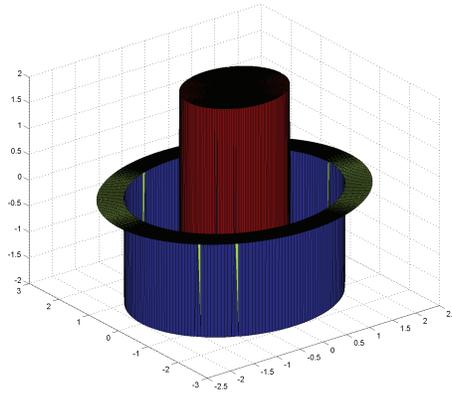}
\caption{Figure of the desired state $y_d$ on the grid of size $h=\frac{14\sqrt{2}}{2^6}$}
\label{example1:yd}
\end{figure}

\begin{figure}[H]
\centering
\subfigure[numerical state $y_{\lambda,h}$ ]{
\includegraphics[width=0.4\textwidth]{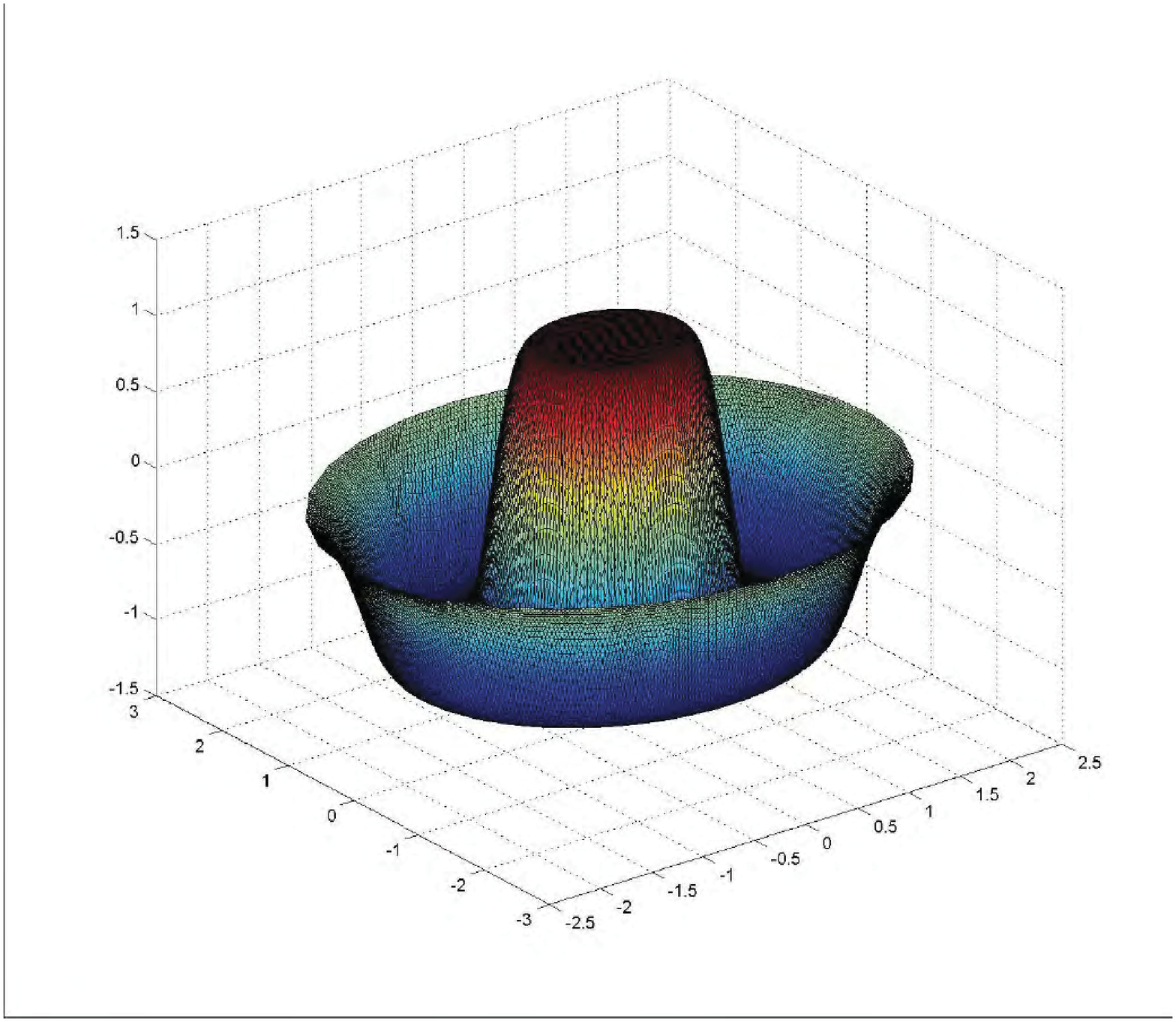}}
\subfigure[numerical control $u_{\lambda,h}$ ]{
\includegraphics[width=0.4\textwidth]{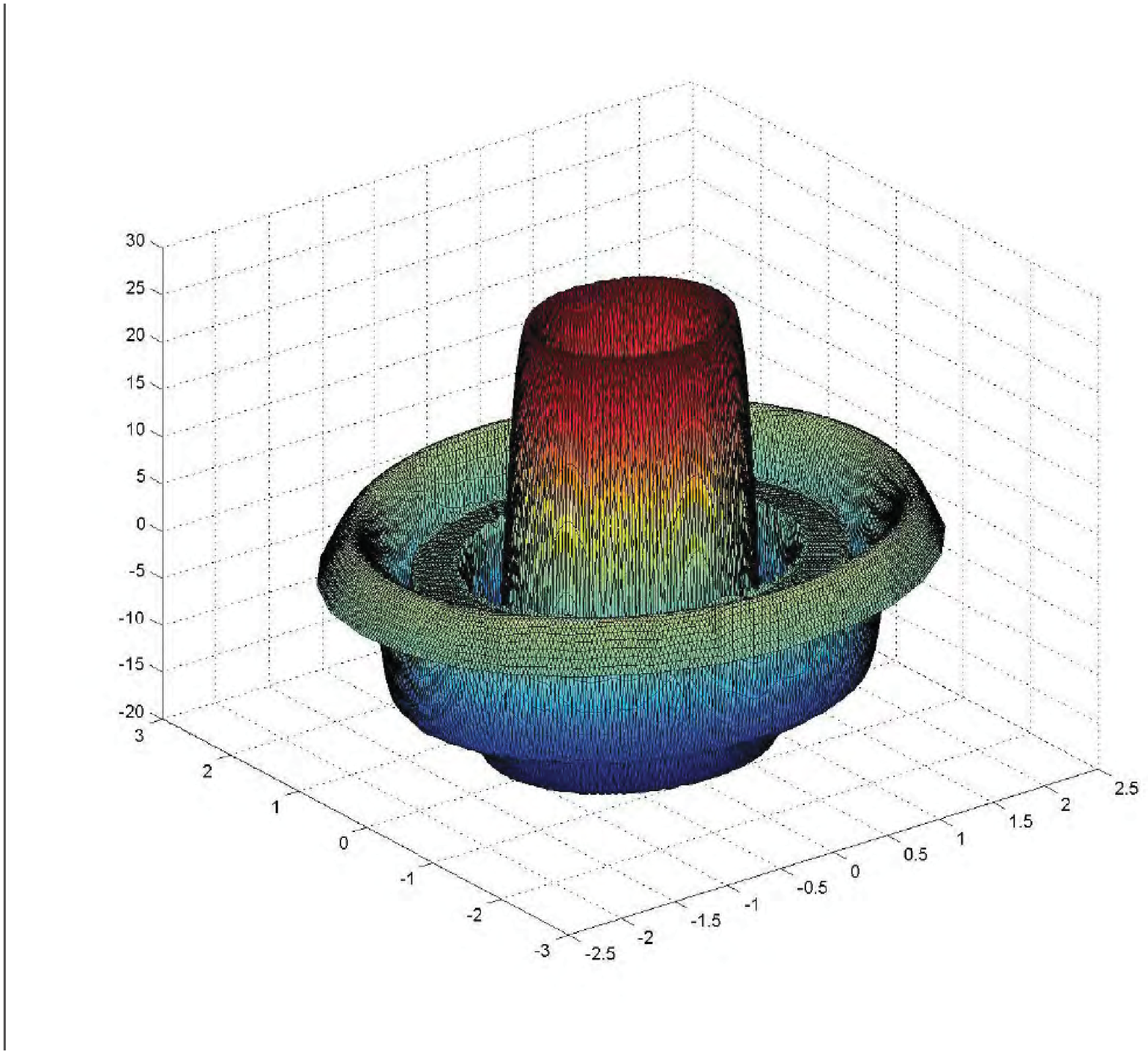}}
\caption{Figures of numerical state and control on the grid of size $h=\frac{14\sqrt{2}}{2^6}$ with $\lambda=10^{-4.5}$}
\label{example1:numerical y and u}
\end{figure}

\begin{table}[H]\small
\ \\
\ \\
\caption{The convergence behavior of our hADMM algorithm and PDAS (a special semi-smooth Newton method) for Example \ref{example:1}.}
\label{table2}
\begin{center}
\begin{tabular}{@{\extracolsep{\fill}}ccccccccccc}
\hline
\multirow{2}{*}{$h$}                        &&\multirow{2}{*}{$\#$dofs}   &&\multirow{2}{*}{$\lambda$}      &&                                  && \multirow{2}{*}{hADMM}           &&\multirow{2}{*}{PDAS}     \\
                                            &&&&&&&&&&\\
\hline
&&&&&&&&&&\\
                                            &&                            &&                                &&\multirow{2}{*}{iter}             &&\multirow{2}{*}{22}                    &&\multirow{2}{*}{17}   \\
                                            &&                            &&\multirow{2}{*}{$10^{-4}$}      &&\multirow{2}{*}{residual $\eta$}        &&\multirow{2}{*}{9.08e-03}          &&\multirow{2}{*}{2.99e-03} \\
                                            &&                            &&                                &&\multirow{2}{*}{time/s}           &&\multirow{2}{*}{1.68}                  &&\multirow{2}{*}{12.20}\\
                                            &&&&&&&&&&\\
                                            &&                            &&                                &&\multirow{2}{*}{iter}             &&\multirow{2}{*}{22}                     &&\multirow{2}{*}{18}  \\
\multirow{2}{*}{$\frac{2.5\sqrt{2}}{2^6}$}  &&\multirow{2}{*}{18977}      &&\multirow{2}{*}{$10^{-4.5}$}    &&\multirow{2}{*}{residual $\eta$}        &&\multirow{2}{*}{9.33e-03}          &&\multirow{2}{*}{3.02e-03} \\
                                            &&                            &&                                &&\multirow{2}{*}{time/s}           &&\multirow{2}{*}{1.65}              &&\multirow{2}{*}{13.06}    \\
                                            &&&&&&&&&&\\
                                            &&                            &&                                &&\multirow{2}{*}{iter}             &&\multirow{2}{*}{22}                    &&\multirow{2}{*}{20}   \\
                                            &&                            &&\multirow{2}{*}{$10^{-5}$}      &&\multirow{2}{*}{residual $\eta$}        &&\multirow{2}{*}{9.49e-03}          &&\multirow{2}{*}{3.04e-03} \\
                                            &&                            &&                                &&\multirow{2}{*}{time/s}           &&\multirow{2}{*}{1.68}                  &&\multirow{2}{*}{14.87}\\
                                            &&&&&&&&&&\\
                                            &&&&&&&&&&\\
\hline
&&&&&&&&&&\\
                                            &&                            &&                                &&\multirow{2}{*}{iter}             &&\multirow{2}{*}{19}                    &&\multirow{2}{*}{34}   \\
                                            &&                            &&\multirow{2}{*}{$10^{-4.5}$}      &&\multirow{2}{*}{residual $\eta$}        &&\multirow{2}{*}{9.02e-03}          &&\multirow{2}{*}{6.13e-03} \\
                                            &&                            &&                                &&\multirow{2}{*}{time/s}           &&\multirow{2}{*}{24.70}                &&\multirow{2}{*}{159.89}\\
                                            &&&&&&&&&&\\
                                            &&                            &&                                &&\multirow{2}{*}{iter}             &&\multirow{2}{*}{19}                    &&\multirow{2}{*}{36}   \\
\multirow{2}{*}{$\frac{2.5\sqrt{2}}{2^7}$}  &&\multirow{2}{*}{76353}      &&\multirow{2}{*}{$10^{-5}$}    &&\multirow{2}{*}{residual $\eta$}        &&\multirow{2}{*}{9.09e-03}          &&\multirow{2}{*}{5.58e-03} \\
                                            &&                            &&                                &&\multirow{2}{*}{time/s}           &&\multirow{2}{*}{24.93}                &&\multirow{2}{*}{171.54}\\
                                            &&&&&&&&&&\\
                                            &&                            &&                                &&\multirow{2}{*}{iter}             &&\multirow{2}{*}{19}                    &&\multirow{2}{*}{37}   \\
                                            &&                            &&\multirow{2}{*}{$10^{-5.5}$}      &&\multirow{2}{*}{residual $\eta$}        &&\multirow{2}{*}{9.12e-03}          &&\multirow{2}{*}{8.35e-03} \\
                                            &&                            &&                                &&\multirow{2}{*}{time/s}           &&\multirow{2}{*}{24.43}                &&\multirow{2}{*}{181.77}\\
                                            &&&&&&&&&&\\
                                            &&&&&&&&&&\\
\hline
&&&&&&&&&&\\
                                            &&                            &&                                &&\multirow{2}{*}{iter}             &&\multirow{2}{*}{19}                    &&\multirow{2}{*}{68}   \\
                                            &&                            &&\multirow{2}{*}{$10^{-5}$}      &&\multirow{2}{*}{residual $\eta$}        &&\multirow{2}{*}{9.03e-03}          &&\multirow{2}{*}{4.67e-03} \\
                                            &&                            &&                                &&\multirow{2}{*}{time/s}           &&\multirow{2}{*}{145.55}               &&\multirow{2}{*}{2938.8}\\
                                            &&&&&&&&&&\\
                                            &&                            &&                                &&\multirow{2}{*}{iter}             &&\multirow{2}{*}{19}                    &&\multirow{2}{*}{71}   \\
\multirow{2}{*}{$\frac{2.5\sqrt{2}}{2^8}$}  &&\multirow{2}{*}{306305}     &&\multirow{2}{*}{$10^{-5.5}$}    &&\multirow{2}{*}{residual $\eta$}        &&\multirow{2}{*}{9.06e-03}          &&\multirow{2}{*}{8.41e-03} \\
                                            &&                            &&                                &&\multirow{2}{*}{time/s}           &&\multirow{2}{*}{145.33}               &&\multirow{2}{*}{3092.9}\\
                                            &&&&&&&&&&\\
                                            &&                            &&                                &&\multirow{2}{*}{iter}             &&\multirow{2}{*}{19}                    &&\multirow{2}{*}{73}   \\
                                            &&                            &&\multirow{2}{*}{$10^{-6}$}      &&\multirow{2}{*}{residual $\eta$}        &&\multirow{2}{*}{9.06e-03}          &&\multirow{2}{*}{4.72e-03} \\
                                            &&                            &&                                &&\multirow{2}{*}{time/s}           &&\multirow{2}{*}{145.57}               &&\multirow{2}{*}{3192.7}\\
                                            &&&&&&&&&&\\
                                            &&&&&&&&&&\\
\hline
\end{tabular}
\end{center}
\end{table}

\begin{table}[H]\small
\ \\
\ \\
\caption{The convergence behavior of our two-phase strategy and PDAS (a special semi-smooth Newton method) for Example \ref{example:1}.}
\label{table3}
\begin{center}
\begin{tabular}{@{\extracolsep{\fill}}ccccccccccc}
\hline
\multirow{2}{*}{$h$}                        &&\multirow{2}{*}{$\#$dofs}   &&\multirow{2}{*}{$\lambda$}      &&                                  && \multirow{2}{*}{two-phase strategy}    &&\multirow{2}{*}{PDAS}\\     &&&&&&&&\multirow{2}{*}{(hADMM $|$ PDAS)}&&\\
                                            &&&&&&&&&&\\
\hline
&&&&&&&&&&\\
                                            &&                            &&                                &&\multirow{2}{*}{iter}             &&\multirow{2}{*}{22 $|$ 5}              &&\multirow{2}{*}{18}   \\
                                            &&                            &&\multirow{2}{*}{$10^{-4}$}      &&\multirow{2}{*}{residual $\eta$}  &&\multirow{2}{*}{9.08e-03 $|$ 4.14e-14}   &&\multirow{2}{*}{4.15e-14} \\
                                            &&                            &&                                &&\multirow{2}{*}{time/s}         &&\multirow{2}{*}{5.33 (1.80 $|$ 3.53)}    &&\multirow{2}{*}{12.92}\\
                                            &&&&&&&&&&\\
                                            &&                            &&                                &&\multirow{2}{*}{iter}             &&\multirow{2}{*}{22 $|$ 6}              &&\multirow{2}{*}{19}   \\
\multirow{2}{*}{$\frac{2.5\sqrt{2}}{2^6}$}  &&\multirow{2}{*}{18977}      &&\multirow{2}{*}{$10^{-4.5}$}    &&\multirow{2}{*}{residual $\eta$}  &&\multirow{2}{*}{9.33e-03 $|$ 3.45e-14}  &&\multirow{2}{*}{3.46e-14} \\
                                            &&                            &&                                &&\multirow{2}{*}{time/s}         &&\multirow{2}{*}{5.81 (1.64 $|$ 4.17)}    &&\multirow{2}{*}{13.74}\\
                                            &&&&&&&&&&\\
                                            &&                            &&                                &&\multirow{2}{*}{iter}             &&\multirow{2}{*}{22 $|$ 8}              &&\multirow{2}{*}{22}   \\
                                            &&                            &&\multirow{2}{*}{$10^{-5}$}      &&\multirow{2}{*}{residual $\eta$}  &&\multirow{2}{*}{9.49e-03 $|$ 4.11e-14}   &&\multirow{2}{*}{4.12e-14} \\
                                            &&                            &&                                &&\multirow{2}{*}{time/s}         &&\multirow{2}{*}{7.32 (1.65 $|$ 5.67)}    &&\multirow{2}{*}{16.35}\\
                                            &&&&&&&&&&\\
                                            &&&&&&&&&&\\
\hline
&&&&&&&&&&\\
                                            &&                            &&                                &&\multirow{2}{*}{iter}             &&\multirow{2}{*}{19 $|$ 8}              &&\multirow{2}{*}{37}  \\
                                            &&                            &&\multirow{2}{*}{$10^{-4.5}$}      &&\multirow{2}{*}{residual $\eta$}  &&\multirow{2}{*}{9.38e-03 $|$ 8.11e-14}   &&\multirow{2}{*}{8.11e-14} \\
                                            &&                            &&                                &&\multirow{2}{*}{time/s}        &&\multirow{2}{*}{62.22 (25.21 $|$ 37.01)}&&\multirow{2}{*}{173.74}\\
                                            &&&&&&&&&&\\
                                            &&                            &&                                &&\multirow{2}{*}{iter}             &&\multirow{2}{*}{19 $|$ 10}             &&\multirow{2}{*}{38}   \\
\multirow{2}{*}{$\frac{2.5\sqrt{2}}{2^7}$}  &&\multirow{2}{*}{76353}      &&\multirow{2}{*}{$10^{-5}$}    &&\multirow{2}{*}{residual $\eta$}  &&\multirow{2}{*}{9.46e-03 $|$ 8.22e-14}   &&\multirow{2}{*}{8.23e-14} \\
                                            &&                            &&                                &&\multirow{2}{*}{time/s}         &&\multirow{2}{*}{69.56 (25.20 $|$ 44.36)}&&\multirow{2}{*}{180.97}\\
                                            &&&&&&&&&&\\
                                            &&                            &&                                &&\multirow{2}{*}{iter}             &&\multirow{2}{*}{19 $|$ 14}             &&\multirow{2}{*}{40}   \\
                                            &&                            &&\multirow{2}{*}{$10^{-5.5}$}      &&\multirow{2}{*}{residual $\eta$}  &&\multirow{2}{*}{9.48e-03 $|$ 8.19e-14}   &&\multirow{2}{*}{8.19e-14} \\
                                            &&                            &&                                &&\multirow{2}{*}{time/s}         &&\multirow{2}{*}{88.30 (25.12 $|$ 63.18)}&&\multirow{2}{*}{195.93}\\
                                            &&&&&&&&&&\\
                                            &&&&&&&&&&\\
\hline
&&&&&&&&&&\\
                                            &&                            &&                                &&\multirow{2}{*}{iter}             &&\multirow{2}{*}{19 $|$ 13}              &&\multirow{2}{*}{73}   \\
                                            &&                            &&\multirow{2}{*}{$10^{-5}$}      &&\multirow{2}{*}{residual $\eta$}  &&\multirow{2}{*}{9.40e-03 $|$ 1.43e-13}   &&\multirow{2}{*}{1.43e-13} \\
                                            &&                            &&                                &&\multirow{2}{*}{time/s}     &&\multirow{2}{*}{679.31 (145.59 $|$ 533.72)}&&\multirow{2}{*}{3152.4}\\
                                            &&&&&&&&&&\\
                                            &&                            &&                                &&\multirow{2}{*}{iter}             &&\multirow{2}{*}{19 $|$ 17}             &&\multirow{2}{*}{75}   \\
\multirow{2}{*}{$\frac{2.5\sqrt{2}}{2^8}$}  &&\multirow{2}{*}{306305}     &&\multirow{2}{*}{$10^{-5.5}$}    &&\multirow{2}{*}{residual $\eta$}  &&\multirow{2}{*}{9.42e-03 $|$ 1.36e-13}   &&\multirow{2}{*}{1.36e-13} \\
                                            &&                            &&                                &&\multirow{2}{*}{time/s}    &&\multirow{2}{*}{835.65 (145.76 $|$ 689.89)}&&\multirow{2}{*}{3241.4}\\
                                            &&&&&&&&&&\\
                                            &&                            &&                                &&\multirow{2}{*}{iter}             &&\multirow{2}{*}{19 $|$ 21}             &&\multirow{2}{*}{77}   \\
                                            &&                            &&\multirow{2}{*}{$10^{-6}$}      &&\multirow{2}{*}{residual $\eta$}  &&\multirow{2}{*}{9.43e-03 $|$ 1.58e-13}   &&\multirow{2}{*}{1.59e-13} \\
                                            &&                            &&                                &&\multirow{2}{*}{time/s}    &&\multirow{2}{*}{986.6 (145.72 $|$ 840.90)}&&\multirow{2}{*}{3334.2}\\
                                            &&&&&&&&&&\\
                                            &&&&&&&&&&\\
\hline
\end{tabular}
\end{center}
\end{table}

\end{xmpl}

\begin{xmpl}\label{example:2}
We consider $\mathrm{\Omega}={[0,14]}^2$ as the test domain and set $\alpha=10^{-3},\ a=-4,\ b=4,\ \sigma=0.5$ and define $g(x)$ as
\begin{equation*}\small
g(x)=\left\{\begin{aligned}
&\frac{1}{6}x^3+\frac{1}{8\pi^3}\cos(2\pi x-\frac{\pi}{2})-\frac{1}{4\pi^2}x\ &x\in \left[0,1\right),\\
&-\left(\frac{1}{6}x^3+\frac{1}{8\pi^3}\cos(2\pi x-\frac{\pi}{2})-x^2+(1-\frac{1}{4\pi^2})x-\frac{1}{3}+\frac{1}{2\pi^2}\right)\ &x\in \left[1,3\right),\\
&\frac{1}{6}x^3+\frac{1}{8\pi^3}\cos(2\pi x-\frac{\pi}{2})-2x^2+(8-\frac{1}{4\pi^2})x-\frac{26}{3}+\frac{1}{\pi^2}\ &x\in \left[3,4\right),\\
&\qquad\qquad\qquad\qquad\qquad\qquad\ 2\ &x\in \left[4,5\right),\\
&-\left(\frac{1}{3}x^3+\frac{1}{4\pi^3}\cos(2\pi x-\frac{\pi}{2})-5x^2+(25-\frac{1}{2\pi^2})x-\frac{131}{3}+\frac{5}{2\pi^2}\right)\ &x\in \left[5,6\right),\\
&\frac{1}{3}x^3+\frac{1}{4\pi^3}\cos(2\pi x-\frac{\pi}{2})-7x^2+(47-\frac{1}{2\pi^2})x-\frac{301}{3}+\frac{7}{2\pi^2}\ &x\in \left[6,8\right),\\
&-\left(\frac{1}{3}x^3+\frac{1}{4\pi^3}\cos(2\pi x-\frac{\pi}{2})-9x^2+(81-\frac{1}{2\pi^2})x-241+\frac{9}{2\pi^2}\right)\ &x\in \left[8,9\right),\\
&\qquad\qquad\qquad\qquad\qquad\qquad-2\ &x\in \left[9,10\right),\\
&\frac{1}{6}x^3+\frac{1}{8\pi^3}\cos(2\pi x-\frac{\pi}{2})-5x^2+(50-\frac{1}{4\pi^2})x-\frac{506}{3}+\frac{5}{2\pi^2}\ &x\in \left[10,11\right),\\
&-\left(\frac{1}{6}x^3+\frac{1}{8\pi^3}\cos(2\pi x-\frac{\pi}{2})-6x^2+(71-\frac{1}{4\pi^2})x-275+\frac{3}{\pi^2}\right)\ &x\in \left[11,13\right),\\
&\frac{1}{6}x^3+\frac{1}{8\pi^3}\cos(2\pi x-\frac{\pi}{2})-7x^2+(98-\frac{1}{4\pi^2})x-\frac{1372}{3}+\frac{7}{2\pi^2}\ &x\in \left[13,14\right].
\end{aligned} \right.
\end{equation*}
Let $y^*(x)=-g(x_1)g(x_2)$,
\begin{equation*}\small
           \mu_a=\left\{ \begin{aligned}
        &0.1\sin(\pi x_1)\sin(\pi x_2)\qquad x\in(4,5)\times(4,5)\ {\rm{or}}\ x\in(9,10)\times(9,10),\\
        &0\qquad\qquad\qquad\qquad\qquad\qquad {\rm{else}},
        \end{aligned} \right.
\end{equation*}
\begin{equation*}\small
           \mu_b=\left\{ \begin{aligned}
        &\hspace{-0.22cm}-0.1\sin(\pi x_1)\sin(\pi x_2)\qquad x\in(4,5)\times(9,10)\ {\rm{or}}\ x\in(9,10)\times(4,5),\\
        &0\qquad\qquad\qquad\qquad\qquad\qquad\qquad {\rm{else}},
        \end{aligned} \right.
\end{equation*}
then from the optimal condition we arrive at
\begin{eqnarray*}
  u^*(x)&=&-\mathrm{\Delta} y^*=g^{(2)}(x_1)g(x_2)+g(x_1)g^{(2)}(x_2), \\
  p\ \ &=&-\alpha u^*, \\
  y_d\ \ &=&y^*+\mu_b-\mu_a+\mathrm{\Delta} p.
\end{eqnarray*}

The exact solution is known in this example and the $L^2$ errors $\|u^*-\overline{u}_{\lambda,h}\|$ on grids of different sizes with nine different values of $\lambda$ from $10^{-2}$ to $10^{-6}$ are given in Table \ref{table4}. As an example, the figures of exact state $y^*$ and numerical state $y_{\lambda,h}$, exact control $u^*$ and numerical control $u_{\lambda,h}$ on the grid of size $h=\frac{14\sqrt{2}}{2^7}$ with $\lambda=10^{-4.5}$ are displayed in Figure \ref{example2:y} and Figure \ref{example2:u}. As stated in Example \ref{example:1}, if a solution with moderate accuracy is sufficient, both hADMM and PDAS are terminated when $\eta_{\mathrm{A}}(\eta_{\mathrm{P}})<10^{-3}$ and the corresponding numerical results are displayed in Table \ref{table5}. Moreover, if more accurate solution is required, we employ the two-phase strategy and compare it with PDAS. Both two algorithms are terminated when $\eta_{\mathrm{A}}(\eta_{\mathrm{P}})<10^{-13}$ in this case and the numerical results are given in Table \ref{table6}.

Table \ref{table4} shows that when $\lambda$ is fixed, the error declines as $h$ decreases until the error is up to a lower bound caused by the regularization. While for a fixed $h$, the error declines as $\lambda$ decreases. The data in Table \ref{table4} verify the error estimates in Section \ref{sec:3}. The last two tables in this example are similar to their counterparts in Example \ref{example:1}. We could find from the numerical results that the hADMM algorithm and the two-phase strategy are faster than PDAS method especially when the finite element grid size $h$ is very small, which verifies the efficiency of the hADMM algorithm and the two-phase strategy.

\begin{table}[H]\footnotesize
\caption{The $L^2$ error $\|u^*-\overline{u}_{\lambda,h}\|$ for Example \ref{example:2}.}
\label{table4}
\begin{center}
\begin{tabular}{@{\extracolsep{\fill}}ccccccccccc}
\hline
\multirow{3}{*}{$\qquad\lambda$} &\multirow{3}{*}{$10^{-2}$} &\multirow{3}{*}{$10^{-2.5}$} &\multirow{3}{*}{$10^{-3}$} &\multirow{3}{*}{$10^{-3.5}$}
&\multirow{3}{*}{$10^{-4}$} &\multirow{3}{*}{$10^{-4.5}$} &\multirow{3}{*}{$10^{-5}$} &\multirow{3}{*}{$10^{-5.5}$} &\multirow{3}{*}{$10^{-6}$}\\
&&&&&&&&&\\
$h$&&&&&&&&&\\
\hline
\multirow{3}{*}{$\frac{14\sqrt{2}}{2^5}$} &\multirow{3}{*}{2.7011}&\multirow{3}{*}{2.6892}&\multirow{3}{*}{2.6935}&\multirow{3}{*}{2.6960}&\multirow{3}{*}{2.6969}&\multirow{3}{*}{2.6973}&\multirow{3}{*}{2.6974}&\multirow{3}{*}{2.6974}&\multirow{3}{*}{2.6975}\\
&&&&&&&&&\\
\multirow{3}{*}{$\frac{14\sqrt{2}}{2^6}$} &\multirow{3}{*}{6.7252e-1}&\multirow{3}{*}{6.4552e-1}&\multirow{3}{*}{6.3471e-1}&\multirow{3}{*}{6.3121e-1}&\multirow{3}{*}{6.3012e-1}&\multirow{3}{*}{6.2978e-1}&\multirow{3}{*}{6.2967e-1}&\multirow{3}{*}{6.2964e-1}&\multirow{3}{*}{6.2963e-1}\\
&&&&&&&&&\\
\multirow{3}{*}{$\frac{14\sqrt{2}}{2^7}$} &\multirow{3}{*}{2.6032e-1}&\multirow{3}{*}{1.5071e-1}&\multirow{3}{*}{1.4671e-1}&\multirow{3}{*}{1.4603e-1}&\multirow{3}{*}{1.4585e-1}&\multirow{3}{*}{1.4580e-1}&\multirow{3}{*}{1.4579e-1}&\multirow{3}{*}{1.4579e-1}&\multirow{3}{*}{1.4578e-1}\\
&&&&&&&&&\\
\multirow{3}{*}{$\frac{14\sqrt{2}}{2^8}$} &\multirow{3}{*}{1.8943e-1}&\multirow{3}{*}{4.5202e-2}&\multirow{3}{*}{3.6722e-2}&\multirow{3}{*}{3.6168e-2}&\multirow{3}{*}{3.6160e-2}&\multirow{3}{*}{3.6139e-2}&\multirow{3}{*}{3.6133e-2}&\multirow{3}{*}{3.6132e-2}&\multirow{3}{*}{3.6131e-2}\\
&&&&&&&&&\\
\multirow{3}{*}{$\frac{14\sqrt{2}}{2^9}$} &\multirow{3}{*}{1.7896e-1}&\multirow{3}{*}{2.2922e-2}&\multirow{3}{*}{9.8562e-3}&\multirow{3}{*}{9.2458e-3}&\multirow{3}{*}{9.1939e-3}&\multirow{3}{*}{9.1874e-3}&\multirow{3}{*}{9.1869e-3}&\multirow{3}{*}{9.1864e-3}&\multirow{3}{*}{9.1863e-3}\\
&&&&&&&&&\\
&&&&&&&&&\\
\hline
\end{tabular}
\end{center}
\end{table}

\begin{figure}[H]
\centering
\subfigure[exact state $y^*$ ]{
\includegraphics[width=0.37\textwidth]{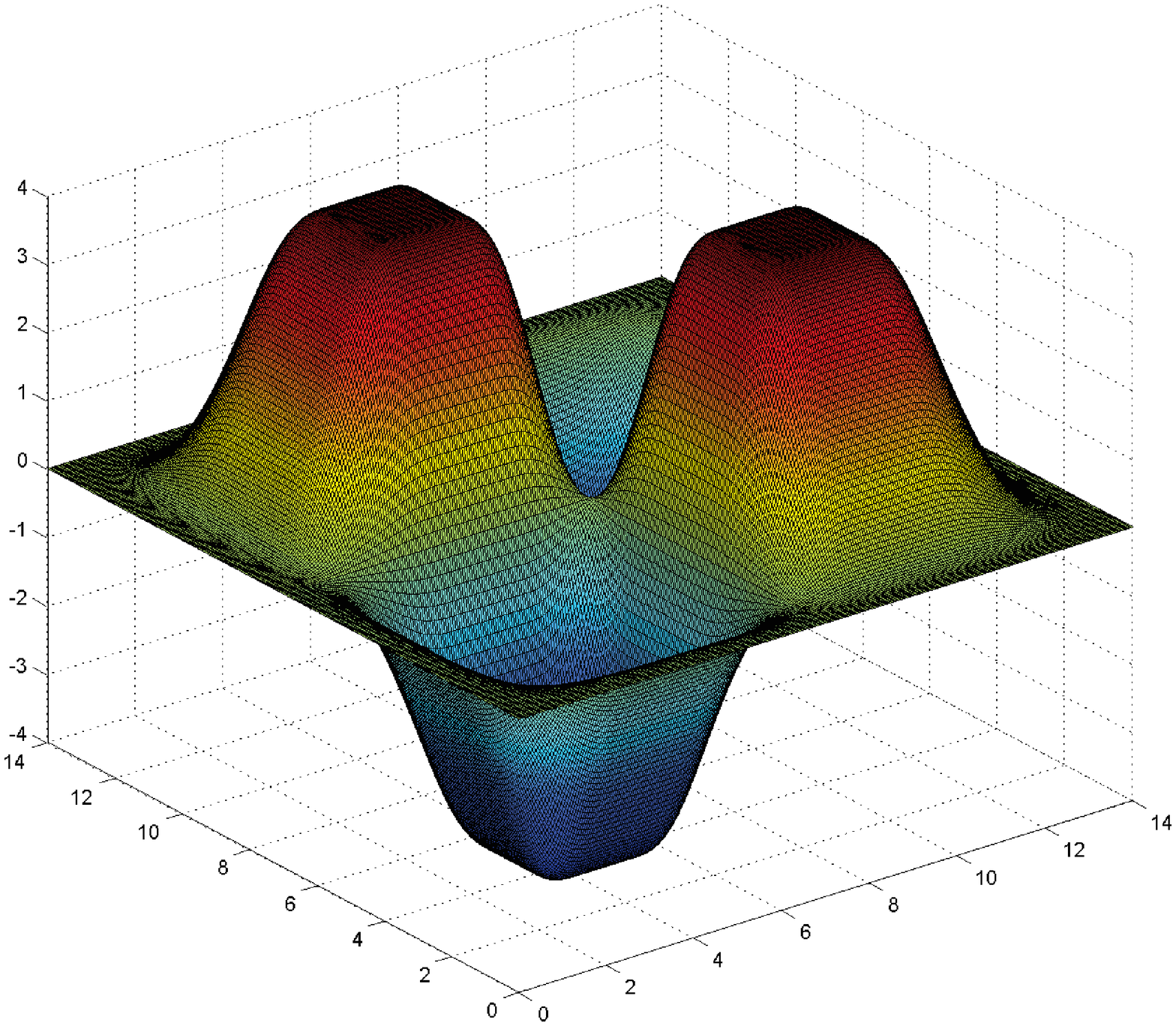}}
\subfigure[numerical state $y_{\lambda,h}$ ]{
\includegraphics[width=0.37\textwidth]{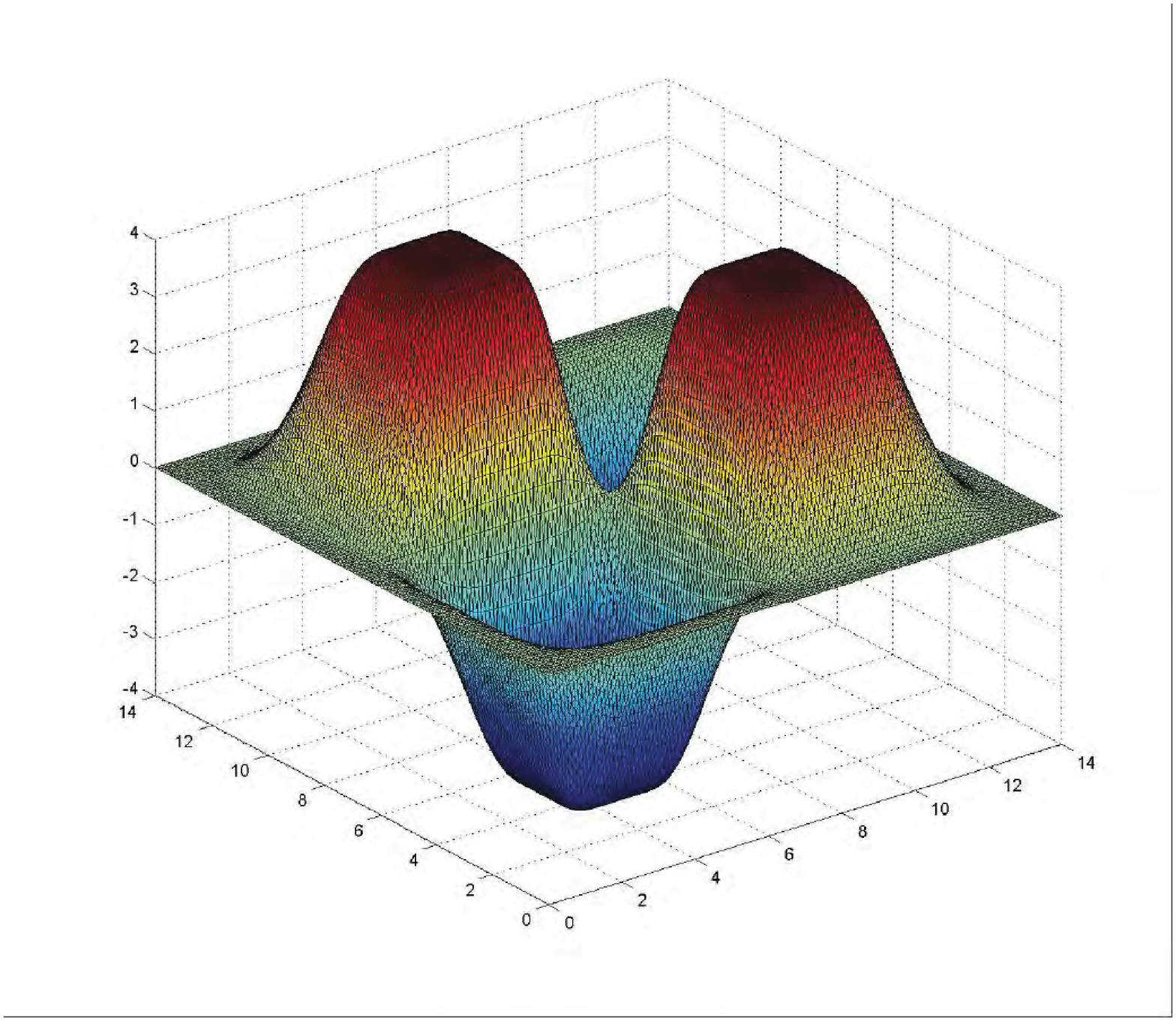}}
\caption{Figures of exact and numerical state on the grid of size $h=\frac{14\sqrt{2}}{2^7}$ with $\lambda=10^{-4.5}$}
\label{example2:y}
\end{figure}

\begin{figure}[H]
\centering
\subfigure[exact control $u^*$ ]{
\includegraphics[width=0.37\textwidth]{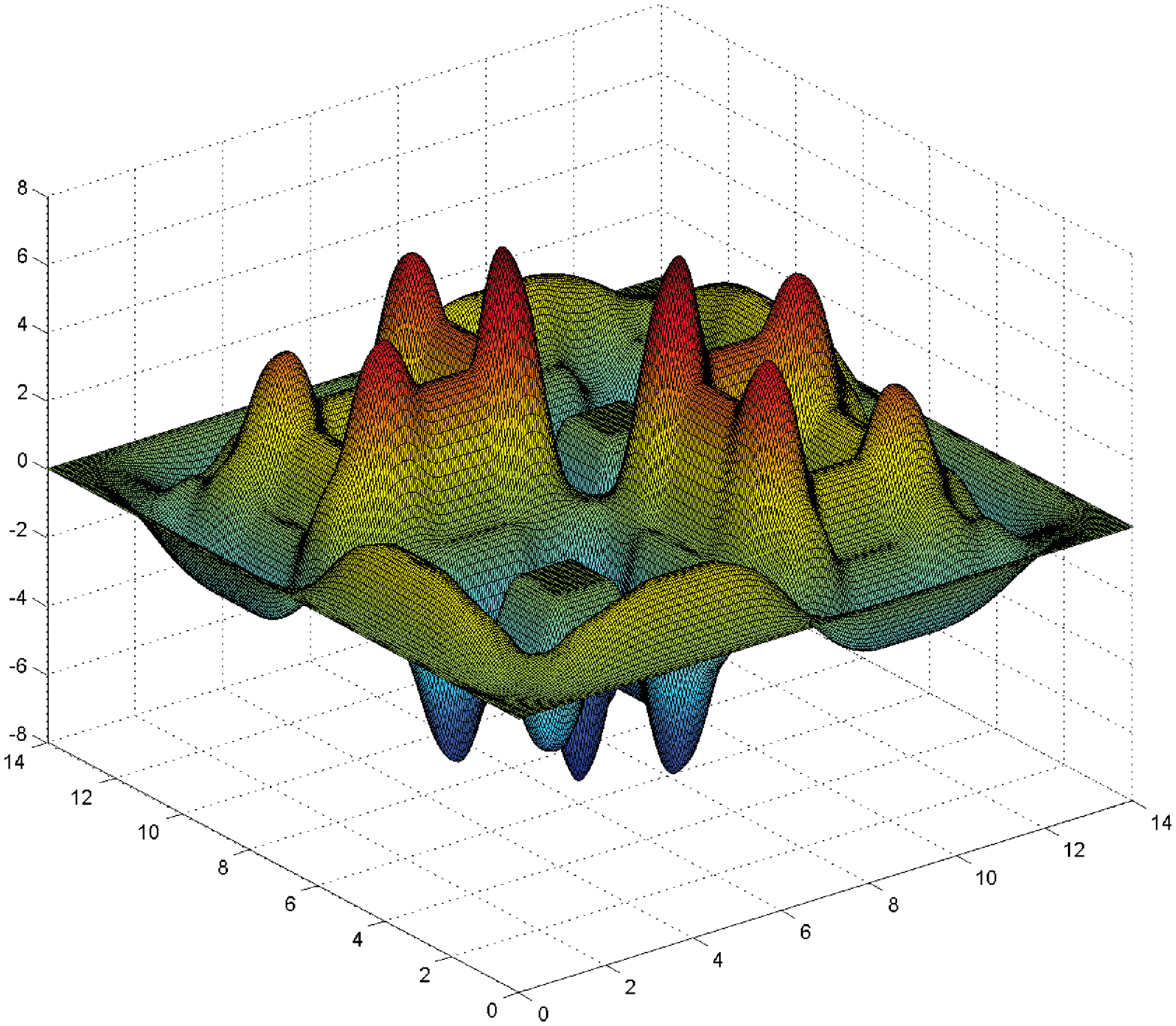}}
\subfigure[numerical control $u_{\lambda,h}$ ]{
\includegraphics[width=0.37\textwidth]{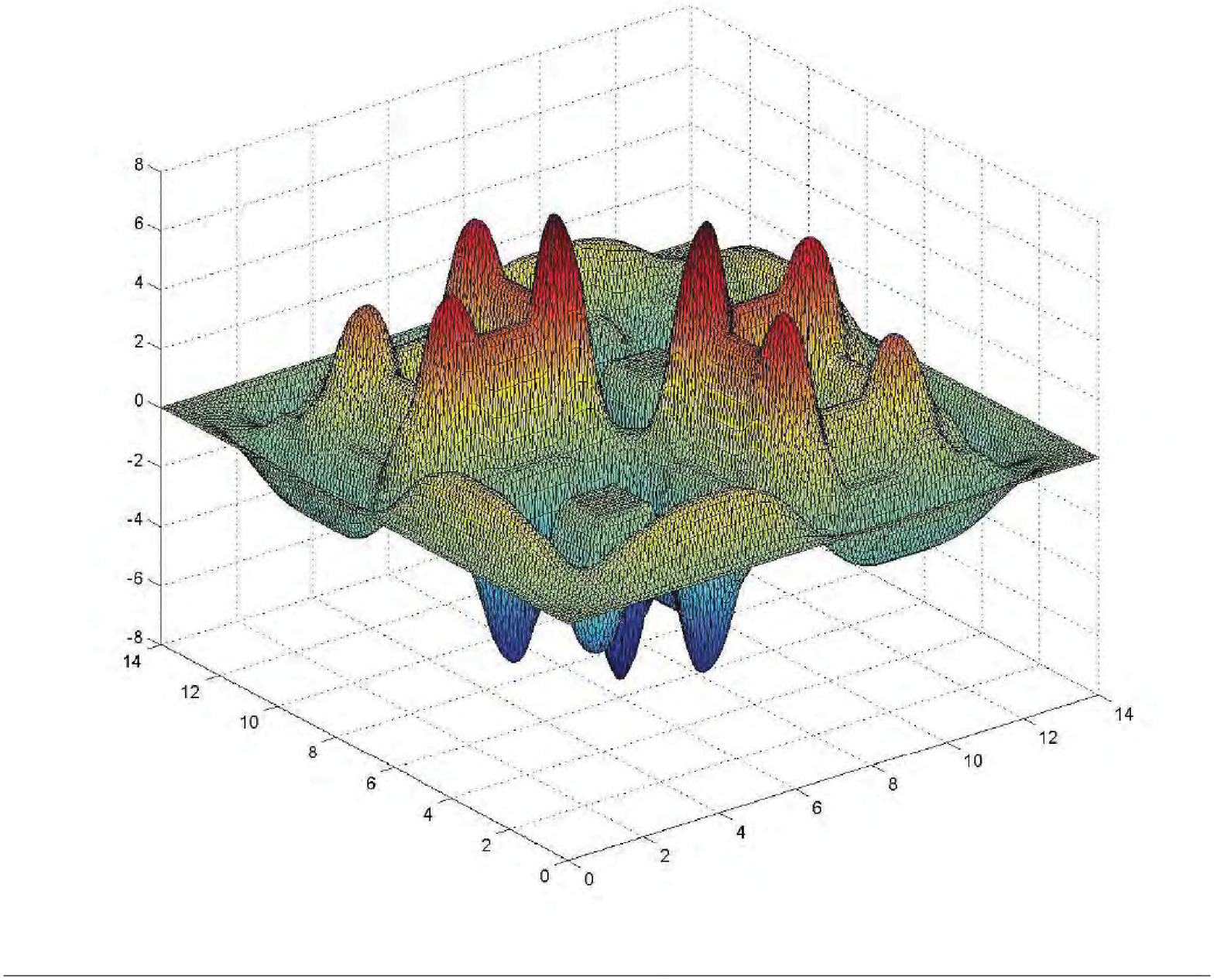}}
\caption{Figures of exact and numerical control on the grid of size $h=\frac{14\sqrt{2}}{2^7}$ with $\lambda=10^{-4.5}$}
\label{example2:u}
\end{figure}

\begin{table}[H]\small
\ \\
\ \\
\caption{The convergence behavior of hADMM and PDAS (a special semi-smooth Newton method) for Example \ref{example:2}.}
\label{table5}
\begin{center}
\begin{tabular}{@{\extracolsep{\fill}}ccccccccccc}
\hline
\multirow{2}{*}{$h$}                        &&\multirow{2}{*}{$\#$dofs}   &&\multirow{2}{*}{$\lambda$}      &&                                  && \multirow{2}{*}{hADMM}               &&\multirow{2}{*}{PDAS} \\
&&&&&&&&&&\\
\hline
&&&&&&&&&&\\
                                            &&                            &&                                &&\multirow{2}{*}{iter}             &&\multirow{2}{*}{31}                    &&\multirow{2}{*}{3}    \\
                                            &&                            &&\multirow{2}{*}{$10^{-4}$}      &&\multirow{2}{*}{residual $\eta$}        &&\multirow{2}{*}{6.84e-04}            &&\multirow{2}{*}{3.66e-04} \\
                                            &&                            &&                                &&\multirow{2}{*}{time/s}           &&\multirow{2}{*}{0.39}                  &&\multirow{2}{*}{1.55} \\
                                            &&&&&&&&&&\\
                                            &&                            &&                                &&\multirow{2}{*}{iter}             &&\multirow{2}{*}{31}                    &&\multirow{2}{*}{3}    \\
\multirow{2}{*}{$\frac{14\sqrt{2}}{2^7}$}   &&\multirow{2}{*}{16129}      &&\multirow{2}{*}{$10^{-4.5}$}    &&\multirow{2}{*}{residual $\eta$}        &&\multirow{2}{*}{7.06e-04}            &&\multirow{2}{*}{4.01e-04} \\
                                            &&                            &&                                &&\multirow{2}{*}{time/s}           &&\multirow{2}{*}{0.36}                  &&\multirow{2}{*}{1.46} \\
                                            &&&&&&&&&&\\
                                            &&                            &&                                &&\multirow{2}{*}{iter}             &&\multirow{2}{*}{30}                    &&\multirow{2}{*}{3}    \\
                                            &&                            &&\multirow{2}{*}{$10^{-5}$}      &&\multirow{2}{*}{residual $\eta$}        &&\multirow{2}{*}{9.12e-04}            &&\multirow{2}{*}{4.12e-04} \\
                                            &&                            &&                                &&\multirow{2}{*}{time/s}           &&\multirow{2}{*}{0.34}                  &&\multirow{2}{*}{1.44} \\
&&&&&&&&&&\\
&&&&&&&&&&\\
\hline
&&&&&&&&&&\\
                                            &&                            &&                                &&\multirow{2}{*}{iter}             &&\multirow{2}{*}{30}                    &&\multirow{2}{*}{5}    \\
                                            &&                            &&\multirow{2}{*}{$10^{-4.5}$}      &&\multirow{2}{*}{residual $\eta$}        &&\multirow{2}{*}{9.39e-04}            &&\multirow{2}{*}{7.81e-04} \\
                                            &&                            &&                                &&\multirow{2}{*}{time/s}           &&\multirow{2}{*}{4.01}                  &&\multirow{2}{*}{14.38}\\
                                            &&&&&&&&&&\\
                                            &&                            &&                                &&\multirow{2}{*}{iter}             &&\multirow{2}{*}{31}                    &&\multirow{2}{*}{5}    \\
\multirow{2}{*}{$\frac{14\sqrt{2}}{2^8}$}   &&\multirow{2}{*}{65025}      &&\multirow{2}{*}{$10^{-5}$}    &&\multirow{2}{*}{residual $\eta$}        &&\multirow{2}{*}{6.48e-04}            &&\multirow{2}{*}{9.17e-04} \\
                                            &&                            &&                                &&\multirow{2}{*}{time/s}           &&\multirow{2}{*}{4.19}                  &&\multirow{2}{*}{14.07}\\
                                            &&&&&&&&&&\\
                                            &&                            &&                                &&\multirow{2}{*}{iter}             &&\multirow{2}{*}{31}                    &&\multirow{2}{*}{5}    \\
                                            &&                            &&\multirow{2}{*}{$10^{-5.5}$}      &&\multirow{2}{*}{residual $\eta$}        &&\multirow{2}{*}{4.69e-04}            &&\multirow{2}{*}{9.58e-04} \\
                                            &&                            &&                                &&\multirow{2}{*}{time/s}           &&\multirow{2}{*}{4.25}                &&\multirow{2}{*}{14.34}\\
&&&&&&&&&&\\
&&&&&&&&&&\\
\hline
&&&&&&&&&&\\
                                            &&                            &&                                &&\multirow{2}{*}{iter}             &&\multirow{2}{*}{30}                    &&\multirow{2}{*}{13}   \\
                                            &&                            &&\multirow{2}{*}{$10^{-5}$}      &&\multirow{2}{*}{residual $\eta$}        &&\multirow{2}{*}{5.76e-04}            &&\multirow{2}{*}{4.66e-04} \\
                                            &&                            &&                                &&\multirow{2}{*}{time/s}           &&\multirow{2}{*}{80.38}                &&\multirow{2}{*}{248.61}\\
&&&&&&&&&&\\
                                            &&                            &&                                &&\multirow{2}{*}{iter}             &&\multirow{2}{*}{30}                    &&\multirow{2}{*}{13}   \\
\multirow{2}{*}{$\frac{14\sqrt{2}}{2^9}$}   &&\multirow{2}{*}{261121}     &&\multirow{2}{*}{$10^{-5.5}$}    &&\multirow{2}{*}{residual $\eta$}        &&\multirow{2}{*}{9.10e-04}            &&\multirow{2}{*}{8.42e-04} \\
                                            &&                            &&                                &&\multirow{2}{*}{time/s}           &&\multirow{2}{*}{80.14}                &&\multirow{2}{*}{245.72}\\
                                            &&&&&&&&&&\\
                                            &&                            &&                                &&\multirow{2}{*}{iter}             &&\multirow{2}{*}{30}                    &&\multirow{2}{*}{13}   \\
                                            &&                            &&\multirow{2}{*}{$10^{-6}$}      &&\multirow{2}{*}{residual $\eta$}        &&\multirow{2}{*}{9.15e-04}            &&\multirow{2}{*}{8.73e-04} \\
                                            &&                            &&                                &&\multirow{2}{*}{time/s}           &&\multirow{2}{*}{81.23}                &&\multirow{2}{*}{247.95}\\
&&&&&&&&&&\\
&&&&&&&&&&\\
\hline
\end{tabular}
\end{center}
\end{table}

\begin{table}[H]\small
\ \\
\ \\
\caption{The convergence behavior of the two-phase strategy and PDAS (a special semi-smooth Newton method) for Example \ref{example:2}.}
\label{table6}
\begin{center}
\begin{tabular}{@{\extracolsep{\fill}}ccccccccccc}
\hline
\multirow{2}{*}{$h$}                        &&\multirow{2}{*}{$\#$dofs}   &&\multirow{2}{*}{$\lambda$}      &&                                  && \multirow{2}{*}{two-phase strategy}   &&\multirow{2}{*}{PDAS}\\
&&&&&&&&\multirow{2}{*}{(hADMM $|$ PDAS)}&&\\
&&&&&&&&&&\\
\hline
&&&&&&&&&&\\
                                            &&                            &&                                &&\multirow{2}{*}{iter}             &&\multirow{2}{*}{31 $|$ 5}              &&\multirow{2}{*}{6}    \\
                                            &&                            &&\multirow{2}{*}{$10^{-4}$}      &&\multirow{2}{*}{residual $\eta$}  &&\multirow{2}{*}{6.84e-04 $|$ 8.84e-14}   &&\multirow{2}{*}{8.67e-14} \\
                                            &&                            &&                                &&\multirow{2}{*}{time/s}         &&\multirow{2}{*}{2.79 (0.36 $|$ 2.43)}    &&\multirow{2}{*}{2.94} \\
                                            &&&&&&&&&&\\
                                            &&                            &&                                &&\multirow{2}{*}{iter}             &&\multirow{2}{*}{31 $|$ 5}              &&\multirow{2}{*}{7}    \\
\multirow{2}{*}{$\frac{14\sqrt{2}}{2^7}$}   &&\multirow{2}{*}{16129}      &&\multirow{2}{*}{$10^{-4.5}$}    &&\multirow{2}{*}{residual $\eta$}  &&\multirow{2}{*}{7.06e-04 $|$ 8.68e-14}   &&\multirow{2}{*}{8.68e-14} \\
                                            &&                            &&                                &&\multirow{2}{*}{time/s}         &&\multirow{2}{*}{2.81 (0.37 $|$ 2.44)}    &&\multirow{2}{*}{3.49} \\
                                            &&&&&&&&&&\\
                                            &&                            &&                                &&\multirow{2}{*}{iter}             &&\multirow{2}{*}{30 $|$ 5}              &&\multirow{2}{*}{6}    \\
                                            &&                            &&\multirow{2}{*}{$10^{-5}$}      &&\multirow{2}{*}{residual $\eta$}  &&\multirow{2}{*}{9.12e-04 $|$ 7.84e-14}   &&\multirow{2}{*}{7.84e-14} \\
                                            &&                            &&                                &&\multirow{2}{*}{time/s}         &&\multirow{2}{*}{2.81 (0.36 $|$ 2.45)}    &&\multirow{2}{*}{2.97} \\

&&&&&&&&&&\\
&&&&&&&&&&\\
\hline
&&&&&&&&&&\\
                                            &&                            &&                                &&\multirow{2}{*}{iter}             &&\multirow{2}{*}{30 $|$ 7}              &&\multirow{2}{*}{12}   \\
                                            &&                            &&\multirow{2}{*}{$10^{-4.5}$}      &&\multirow{2}{*}{residual $\eta$}  &&\multirow{2}{*}{9.39e-04 $|$ 1.74e-13}   &&\multirow{2}{*}{1.75e-13} \\
                                            &&                            &&                                &&\multirow{2}{*}{time/s}         &&\multirow{2}{*}{25.33 (4.81 $|$ 20.52)}  &&\multirow{2}{*}{35.86}\\
                                            &&&&&&&&&&\\
                                            &&                            &&                                &&\multirow{2}{*}{iter}             &&\multirow{2}{*}{31 $|$ 7}              &&\multirow{2}{*}{11}   \\
\multirow{2}{*}{$\frac{14\sqrt{2}}{2^8}$}   &&\multirow{2}{*}{65025}      &&\multirow{2}{*}{$10^{-5}$}    &&\multirow{2}{*}{residual $\eta$}  &&\multirow{2}{*}{6.48e-04 $|$ 1.73e-13}   &&\multirow{2}{*}{1.74e-13} \\
                                            &&                            &&                                &&\multirow{2}{*}{time/s}         &&\multirow{2}{*}{25.12 (4.53 $|$ 20.59)}  &&\multirow{2}{*}{32.46}\\
                                            &&&&&&&&&&\\
                                            &&                            &&                                &&\multirow{2}{*}{iter}             &&\multirow{2}{*}{31 $|$ 7}              &&\multirow{2}{*}{11}   \\
                                            &&                            &&\multirow{2}{*}{$10^{-5.5}$}      &&\multirow{2}{*}{residual $\eta$}  &&\multirow{2}{*}{4.69e-04 $|$ 1.73e-13}   &&\multirow{2}{*}{1.73e-13} \\
                                            &&                            &&                                &&\multirow{2}{*}{time/s}         &&\multirow{2}{*}{25.45 (4.78 $|$ 20.67)}  &&\multirow{2}{*}{32.70}\\
&&&&&&&&&&\\
&&&&&&&&&&\\
\hline
&&&&&&&&&&\\
                                            &&                            &&                                &&\multirow{2}{*}{iter}             &&\multirow{2}{*}{30 $|$ 14}              &&\multirow{2}{*}{24}  \\
                                            &&                            &&\multirow{2}{*}{$10^{-5}$}      &&\multirow{2}{*}{residual $\eta$}  &&\multirow{2}{*}{5.76e-04 $|$ 3.37e-13}   &&\multirow{2}{*}{3.37e-13} \\
                                            &&                            &&                                &&\multirow{2}{*}{time/s}      &&\multirow{2}{*}{344.38 (80.94 $|$ 263.44)}&&\multirow{2}{*}{462.94}\\
                                            &&&&&&&&&&\\
                                            &&                            &&                                &&\multirow{2}{*}{iter}             &&\multirow{2}{*}{30 $|$ 14}             &&\multirow{2}{*}{25}   \\
\multirow{2}{*}{$\frac{14\sqrt{2}}{2^9}$}   &&\multirow{2}{*}{261121}     &&\multirow{2}{*}{$10^{-5.5}$}    &&\multirow{2}{*}{residual $\eta$}  &&\multirow{2}{*}{9.10e-04 $|$ 3.43e-13}   &&\multirow{2}{*}{3.43e-13} \\
                                            &&                            &&                                &&\multirow{2}{*}{time/s}      &&\multirow{2}{*}{345.09 (80.94 $|$ 264.15)}&&\multirow{2}{*}{479.71}\\
                                            &&&&&&&&&&\\
                                            &&                            &&                                &&\multirow{2}{*}{iter}             &&\multirow{2}{*}{30 $|$ 14}             &&\multirow{2}{*}{26}   \\
                                            &&                            &&\multirow{2}{*}{$10^{-6}$}      &&\multirow{2}{*}{residual $\eta$}  &&\multirow{2}{*}{9.15e-04 $|$ 3.35e-13}   &&\multirow{2}{*}{3.35e-13} \\
                                            &&                            &&                                &&\multirow{2}{*}{time/s}     &&\multirow{2}{*}{346.69 (80.28 $|$ 266.41)}&&\multirow{2}{*}{496.61}\\
&&&&&&&&&&\\
&&&&&&&&&&\\
\hline
\end{tabular}
\end{center}
\end{table}

\end{xmpl}

\section{Conclusion}
\label{sec:6}
In this paper, state-constrained elliptic control problems are considered, where the Lagrange multipliers associated to the state constraints are only measure functions. To tackle this difficulty, we utilize Lavrentiev regularization. After that, the regularized problem is discretized by full finite element discretization, in which both the state and control are discretized by piecewise linear functions. We derive error analysis of the overall error resulted from regularization and discretization. To solve the discretized problem efficiently, a heterogeneous alternating direction method of multipliers (hADMM) is proposed. If more accurate solution is required, a two-phase strategy is proposed, in which the primal-dual active set (PDAS) method is used as a postprocessor of the hADMM. Numerical results not only verify the analysis results of error estimate but also show the efficiency of the proposed algorithm.\\
\\
\\
\emph{Acknowledgments}. We would like to thank Dr. Long Chen very much for the contribution of the FEM package iFEM \cite{Chen} in Matlab. Also we are grateful for valuable suggestions of our colleagues.


\end{document}